\documentclass[english,12pt,pre,unsortedaddress,nofootinbib,superscriptaddress]{revtex4-1}
\usepackage{ae,aecompl}
\usepackage[T1]{fontenc}
\usepackage[latin9]{inputenc}
\usepackage{float}
\usepackage{amsmath,amssymb}
\usepackage{mathtools}
\usepackage{array}
\usepackage{graphicx}
\usepackage{color}
\usepackage{subcaption}
\usepackage{ragged2e}
\DeclareCaptionJustification{justified}{\justifying}
\usepackage{enumitem}
\usepackage{setspace}

\makeatletter

\floatstyle{ruled}
\newfloat{algorithm}{H}{loa}
\providecommand{\algorithmname}{Algorithm}
\floatname{algorithm}{\protect\algorithmname}

\makeatother

\usepackage{babel}

\begin{document}

\newcommand{\norm}[1]{\left\lVert#1\right\rVert}
\global\long\def\bm#1{\boldsymbol{#1}}
\global\long\def\abs#1{\left|#1\right|}
\newcommand{\V}[1]{\boldsymbol{#1}} 
\global\long\def\grad{\V{\nabla}}
\newcommand{\D}[1]{\Delta #1}
\global\long\def\rnf{r_\text{nf}}
\global\long\def\kpar{k_{\parallel}}
\global\long\def\kparv{\V{k}_{\parallel}}
\global\long\def\rcut{r_{\text{cut}}}
\global\long\def\four#1{\widehat{#1}}
\global\long\def\olap#1{{#1}_\text{over}}
\global\long\def\inside#1{{#1}_\text{in}}
\global\long\def\near#1{#1^{(n)}}
\global\long\def\far#1{#1^{(f)}}
\global\long\def\GBC{G_\text{BC}}
\global\long\def\GDP{G_\text{DP}}
\global\long\def\fimg{f^\text{(img)}}
\global\long\def\farf{\widetilde{f}}
\global\long\def\img#1{{#1}^\text{(img)}}
\global\long\def\olap#1{{#1}_{<2H_E}}
\global\long\def\nolap#1{{#1}_{>2H_E}}
\global\long\def\cor#1{{#1}^{(c)}}
\global\long\def\kT{k_BT}
\newcommand\ddfrac[2]{\frac{\displaystyle #1}{\displaystyle #2}}
\newcommand{\angstrom}{\text{\normalfont\AA}}
\newcommand{\widebar}[1]{%
   \hbox{%
     \vbox{%
       \hrule height 0.5pt 
       \kern0.5ex
       \hbox{%
         \kern-0.1em
         \ensuremath{#1}%
         \kern-0.1em
       }%
     }%
   }%
}

\newcommand{\Donev}[1]{{\bf[{\color{red}#1}]}}
\newcommand{\Raul}[1]{{\bf[{\color{blue}#1}]}}
\newcommand{\deleted}[1]{}
\newcommand{\red}[1]{{\color{red}#1}}
\newcommand{\rev}[1]{{\color{black}#1}}

\title{A fast spectral method for electrostatics in doubly-periodic slit channels}
\author{Ondrej Maxian}
\affiliation{Courant Institute, NYU, New York, NY 10012}
\author{Raul P. Pel\'aez}
\affiliation{Department of Theoretical Condensed Matter Physics, Universidad Aut\'onoma de Madrid, 28049, Madrid, Spain}
\author{Leslie Greengard}
\affiliation{Courant Institute, NYU, New York, NY 10012}
\affiliation{Center for Computational Mathematics, Flatiron Institute, New York, NY 10010}
\author{ Aleksandar Donev}
\affiliation{Courant Institute, NYU, New York, NY 10012}

\begin{abstract}
We develop a fast method for computing the electrostatic energy and forces for a collection of charges in doubly-periodic slabs with jumps in the dielectric permittivity at the slab boundaries. Our method achieves spectral accuracy by using Ewald splitting to replace the original Poisson equation for nearly-singular sources with a smooth far-field Poisson equation, combined with a localized near-field correction. Unlike existing spectral Ewald methods, which make use of the Fourier transform in the aperiodic direction, we recast the problem as a two-point boundary value problem in the aperiodic direction for each transverse Fourier mode, for which exact analytic boundary conditions are available. We solve each of these boundary value problems using a fast, well-conditioned Chebyshev method. In the presence of dielectric jumps, combining Ewald splitting with the classical method of images results in smoothed charge distributions which overlap the dielectric boundaries themselves. We show how to preserve \rev{spectral} accuracy in this case through the use of a harmonic correction which involves solving a simple Laplace equation with smooth boundary data. We implement our method on Graphical Processing Units, and combine our doubly-periodic Poisson solver with Brownian Dynamics to study the equilibrium structure of double layers in binary electrolytes confined by dielectric boundaries. Consistent with prior studies, we find strong charge depletion near the interfaces due to repulsive interactions with image charges, which points to the need for incorporating polarization effects in understanding confined electrolytes, both theoretically and computationally. 
\end{abstract}

\maketitle

\section{Introduction}
The evaluation of electrostatic interactions in a collection of charges is a classical problem in computational physics, with applications to the study of electrolyte solutions, macromolecules, ion channels, and other systems. In Molecular Dynamics (MD) and Brownian Dynamics (BD) methods, the forces between the charges need to be computed at least once per time step, while in Monte Carlo (MC) methods the total electrostatic energy of the entire collection of charges is required. Hybrid BD-MC schemes require computing both forces and energy. Each of the many linear-scaling methods available to compute these quantities (for varying geometries and boundary conditions) falls into one of two categories: fast multipole methods \cite{FMM_Review_Greengard,FMM_MD_Tornberg,liang2020harmonic}, and variants of particle-(particle-particle)-mesh (P3M) methods, including pre-corrected FFT and spectral Ewald (SE) methods \cite{DoublyPeriodic_Poisson,SpectralEwald_NFFT,SpectralEwald_Electrostatics,SpectralEwaldElectrostatics_Review,DISCOS_Periodic,precorrectedFFT}. Our focus in this work is on the latter type of method, which tends to be more efficient than the former for homogeneous charge distributions. 

In this paper we develop a linear-scaling spectral Ewald method to compute electrostatic forces and energy in \emph{doubly-periodic} or \emph{slab} geometries with dielectric jumps at the boundaries of the slab, which is applicable to the study of (Debye) double layers in confined electrolyte solutions. Specifically, we will consider charges in a slab $0<z<H$ immersed in a uniform dielectric medium in a domain that is periodic in the $xy$ plane but potentially has jumps in the dielectric permittivity at $z=0$ and $z=H$. We will consider not point charges, but rather Gaussian charges of width $g_w>0$. This approach avoids divergent self-interactions and unnecessarily stiff electrostatic interactions at short distances and is consistent with our focus on electrolyte solutions, for which the charged particles are solvated ions that are not actually point charges to begin with. Nevertheless, our method applies to arbitrarily small $g_w$ (to within roundoff errors), allowing us to approach the limit of point charges if desired.

There are already a large number of methods in the literature for doubly-periodic electrostatics \cite{DPPoisson_ELC,DoublyPeriodic_Poisson,DielectricDoublyPeriodicSlabs,DPPoisson_Ewald3D,DPPoisson_MMM2D,SpectralEwald_NFFT,DoublyPeriodic_Electrostatics,liang2020harmonic}. We will not attempt to review and compare all of them, but focus instead on highlighting the improvements of our approach over recent SE methods \cite{SpectralEwald_NFFT,DoublyPeriodic_Electrostatics}, to which our approach is most closely related. Our method combines a number of ideas from the existing literature with some new ideas and new numerical methods, which results in an improved P3M method in the slab geometry. Like many previous approaches, we use the method of images to tackle the presence of dielectric jumps \cite{DPPoisson_ELC,DPPoisson_MMM2D,liang2020harmonic} and convert the problem to a doubly-periodic problem in a medium with uniform dielectric permittivity. At the same time, however, our method avoids the need to consider the full (infinite) image system by imposing exact boundary conditions at $z=0$ and $z=H$, and using a grid-based Poisson/Laplace solver to account for the distant images. We make key use of Ewald splitting, so that some aspects of our method are very similar to existing SE approaches.

Our use of Ewald splitting is different from that in \cite{SpectralEwald_NFFT,DoublyPeriodic_Electrostatics}.  The traditional view is to consider Ewald splitting as an analytical technique that separates the electrostatic sums for 
point-like charges into ``real-space" and ``Fourier-space" contributions. 
The real-space or \emph{near-field} part is easy to handle by direct summation over pairs of nearby particles, while the Fourier-space or \emph{far-field} part is handled using Fast Fourier Transform (FFT) based methods. These are
simplest to understand in triply-periodic domains \cite{SpectralEwald_Electrostatics}, and the corresponding formulas are easy to derive.
This standard view  has been extended to other geometries including slabs, as reviewed in \cite{Ewald_Electrostatics_Review,SpectralEwaldElectrostatics_Review}.

One drawback of the traditional view of Ewald splitting is that, by focusing on the Fourier domain alone, non-physical sampling requirements are sometimes imposed on the method.
The paper \cite{DoublyPeriodic_Electrostatics}, which was perhaps the 
first to introduce and carefully analyze a spectrally accurate approach for 
doubly-periodic geometries, requires large oversampling 
factors that depend on both the tolerance and the aspect ratio of the domain.
While a later paper \cite{SpectralEwald_NFFT} was able to overcome this,
we believe the applicability of Ewald-type methods to more general geometries is easier to understand in the framework we present here, which makes simultaneous use of both Fourier and more standard PDE-based ideas.
 This framework is key to understanding the novel parts of our algorithm. 

Thus, we begin by considering Ewald splitting as an improvement of P3M approaches, where the first step is to create a smooth source distribution by convolution with a Gaussian, the second step involves solving the Poisson equation with this smooth right-hand side, and the third step is to correct for the (localized) errors introduced by the initial smoothing. One can use any sufficiently accurate solver for the smoothed problem, not necessarily one based on Fourier analysis.
The third step is handled (as in any Ewald-type method)
analytically, using pairwise summation over some collection of near neighbors.

Having split the problem into a PDE with smooth data (accurately representing 
the far-field interactions) and a near-field correction, it remains only to correctly specify appropriate boundary conditions in the unbounded $z$ direction and to construct an appropriate fast solver. It is straightforward to derive a Robin-type condition on the slab boundary for each mode, through an analysis of the Dirichlet-to-Neumann (DtN) map.\footnote{The use of the DtN map is not new; see, for example, Eqs.\ (17,18) in \cite{DoublyPeriodic_Poisson} or the FMM-based algorithm in \cite{liang2020harmonic}, but our use of it within the SE framework appears not to have been explored.} After Fourier transformation with respect to the periodic directions, we use a real-space Chebyshev spectral solver in the now finite $z$ direction. Unlike pure Fourier-based schemes, such as
\cite{DoublyPeriodic_Electrostatics}, all that is required is that the source
distribution be {\em resolved} by the Chebyshev grid. 
No oversampling is required and the aspect ratio of the domain plays no role.
Unlike \cite{SpectralEwald_NFFT},
we maintain spectral accuracy rather than switching over to an algebraic convergence rate. That said, since Chebyshev methods are easily implemented efficiently and robustly using a Fourier transform on a double-sized grid, we still make use of three-dimensional FFTs as the key component to achieving $O(N \log N)$ scaling where $N$ is the number of unknowns. The \emph{fixed} oversampling factor of 2 is typical for any FFT-based aperiodic convolution. 

For a single dielectric interface, our Ewald splitting approach can be applied \emph{after} using the classical image construction to handle the dielectric jump in the $z$ direction. In particular, a set of images can be constructed and Ewald splitting can be used to smear the charges, thereby allowing for a coarse grid in the doubly-periodic Poisson solver. This approach runs into difficulties, however, for slabs with two (or more) dielectric jumps. Our method avoids the inherent problem of having infinitely many images for a slab geometry by only including the images of charges that are sufficiently close to the dielectric boundary for the corresponding Ewald cloud to overlap the boundary. We account for the rest of the infinite image system through a harmonic correction potential that can be computed analytically from the mismatch in the boundary conditions at the slab boundaries. This simple but powerful idea appears have been overlooked in the field, and is easily combined with other doubly-periodic Poisson solvers, such as the nonuniform FFT-based solver proposed in \cite{SpectralEwald_NFFT}.

An outline of the paper follows. In Section \ref{sec:problem}, we present the mathematical formulation of the problem. In Section \ref{sec:smoothdp}, we develop a novel continuum approach for smooth doubly-periodic electrostatics problems based on the DtN map. We then present our variant of Ewald splitting in Section \ref{sec:ewnoBC} that maps a nonsmooth problem into a smooth one. In Section \ref{sec:ourstuff}, we present our main contribution: a continuum approach to doubly-periodic electrostatics for slabs that combines a restricted image construction, Ewald splitting, the DtN-based doubly-periodic solver, and a correction approach for coarse-graining images that are sufficiently far away from the slab. In Section \ref{sec:numerics}, we present a discretization of the continuum formulation using FFTs in the $xy$ plane and Chebyshev polynomials in the $z$ direction, which we implement in a public-domain code running on Graphical Processing Units (GPUs). In Section \ref{sec:tests}, we validate the accuracy of our method by comparing to reference analytical and numerical results, and in Section \ref{sec:raul} we combine our electrostatic solver with Brownian Dynamics to study binary electrolyte solutions in slit channels with either uncharged (Section \ref{sec:uncharged}) or charged walls (Section \ref{sec:charged}). By comparing our results to reference Monte Carlo results from the literature and analytical solutions, we validate our method and establish the importance of polarization effects that come from jumps in dielectric permittivity. We conclude with a summary and a discussion of future directions in Section \ref{sec:conclusions}.

\section{Problem statement \label{sec:problem}}
We consider solving an electrostatics problem for a collection of $N$ Gaussian charges with strengths $q_k$ and positions $\V{z}_k$. The domain geometry is that of a slit channel: periodic in the $x$ and $y$ directions on $[0,L_x]$ and $[0,L_y]$, respectively, and unbounded in $z$. We will also assume that the charges are contained within a finite region $z \in (0,H)$, and that there are fixed surface charge densities $\sigma_{b/t}(x,y)$ on the bottom and top boundaries of this region ($z=0$ and $z=H$). These assumptions give the electrostatic equation for the potential $\phi(\V{x}=(x,y,z))$,
\begin{gather}
\label{eq:first}
    - \grad \cdot \left(\widetilde{\epsilon} \grad \phi \right) = \rho, \quad \text{where}\\
    \label{eq:rho}
    \rho(\V{x}) =\rho(x,y,z) = \sum_{k=1}^N q_k g\left(\norm{\V{x}-\V{z}_k}\right) + \sigma_b(x,y) \delta(z) + \sigma_t(x,y)\delta(z-H), 
\end{gather}
on $(x,y) \in [0,L_x] \times [0,L_y]$ with $z$ unbounded. We assume that each charge has a Gaussian charge density
\begin{equation}
     \label{eq:gcharge}
g(r) = \frac{1}{\left(2\pi g_w^2 \right)^{3/2}} \exp{\left(-\frac{r^2}{2 g_w^2}\right)},
\end{equation}
with standard deviation $g_w$ related to the physical size of the charges, with $g_w \rightarrow 0$ corresponding to point charges. In\ \eqref{eq:first}, the potential $\phi$ is only unique up to a constant, and so we set 
\begin{equation}
    \label{eq:phi00}
    \phi(\V{0})=0.
\end{equation}

Our goal is to solve\ \eqref{eq:first} for a slab with piecewise constant dielectric permittivity in the $z$ direction,
\begin{equation}
    \widetilde{\epsilon} (\V{x}) = \widetilde{\epsilon}(z) = 
    \begin{cases} 
    \epsilon_b & z < 0 \\ 
    \epsilon & 0 < z < H \\ 
    \epsilon_t & z > H 
    \end{cases}
\end{equation}
Substituting this into the electrostatic equation\ \eqref{eq:first}, we obtain a Poisson equation for $z \in (0,H)$, 
\begin{gather}
\label{eq:poiss}
\epsilon \Delta \phi(\V{x}) = -\sum_{k=1}^N q_k g\left(\norm{\V{x}-\V{z}_k}\right):=-f(\V{x})
\end{gather}
together with boundary conditions on the potential and electric displacement at $z=0$ and $z=H$, 
\begin{gather}
\label{eq:cbc1}
\phi(x,y,z\rightarrow 0^+)=\phi(x,y,z\rightarrow 0^-)\\[2 pt]
\label{eq:cbc2}
\phi(x,y,z\rightarrow H^-)=\phi(x,y,z\rightarrow H^+)\\[2 pt]
\label{eq:cbc3}
\epsilon \frac{\partial \phi}{\partial z}(x,y,z\rightarrow 0^+)-\epsilon_b \frac{\partial \phi}{\partial z}(x,y,z\rightarrow 0^-)+\sigma_b(x,y)=0,\\[2 pt]
\label{eq:cbc4}
\epsilon \frac{\partial \phi}{\partial z}(x,y,z\rightarrow H^-)-\epsilon_t \frac{\partial \phi}{\partial z}(x,y,z\rightarrow H^+)-\sigma_t(x,y)=0. 
\end{gather}

We assume here that the domain is overall electroneutral (including the wall-bound charge densities), 
\begin{equation}
    \label{eq:elecneu}
    \sum_{k=1}^N q_k +\int_{0}^{L_x} \int_{0}^{L_y}  \left(\sigma_b(x,y)+\sigma_t(x,y)\right) \, dx \, dy  = 0, 
\end{equation}
and so the electric field $\grad \phi$ must decay to zero as $z \rightarrow \pm \infty$. In\ \eqref{eq:elecneu}, we have assumed that each Gaussian charge density\ \eqref{eq:gcharge} is fully contained inside the slab, so that it integrates to unity on $[0,H]$. Of course, this is not possible exactly since a Gaussian is not compactly supported. To address this, we truncate the Gaussian at a finite distance $n_\sigma g_w \geq 4g_w$, so that the three-dimensional integral of each charge density is at least $99.9\%$ of the charge $q$. For the Gaussian to be fully contained inside the slab, we assume that the truncated Gaussian envelopes do not overlap the dielectric boundaries, $z_k \in [n_\sigma g_w, H-n_\sigma g_w]$.

For MCMC simulations, we need to compute the electrostatic energy 
\begin{equation}
    U = \frac{1}{2} \int_{0}^{L_x} \int_{0}^{L_y} \int_{-\infty}^\infty \phi(x,y,z) \rho(x,y,z) \, dz \, dy \, dx.  
\end{equation}
Substituting the expression\ \eqref{eq:rho} for the charge density $\rho$, we obtain the energy
\begin{gather}
    \label{eq:Ufinal}
U = \frac{1}{2} \sum_{k=1}^N q_k \widebar{\phi}(\bm{z}_k) + \frac{1}{2}\int_{0}^{L_y} \int_{0}^{L_x} (\sigma(x,y,0)\phi(x,y,0)+\sigma(x,y,H)\phi(x,y,H))\, dx \, dy,\\
\text{where} \quad
\label{eq:pbar}
\widebar{\phi}(\bm{z}) = \int_{\bm{x}} \phi(\bm{x}) g\left(\norm{\V{x}-\V{z}}\right) \, d\bm{x}
\end{gather}
is the convolution of the pointwise potential $\phi(\V{x})$ with a Gaussian of width $g_w$. This means that the energy\ \eqref{eq:Ufinal} is the sum of the average potential at the center of each charge plus the energy due to the surface charge densities $\sigma_{b/t}$. 

The electrostatic forces $\V{F}_k=-\partial U/\partial \V{z}_k$ can now be determined in a straightforward way from the energy\ \eqref{eq:Ufinal}. The pointwise electric field can be determined from the pointwise potential by
\begin{equation}
\bm{E} = - \grad \phi,
\end{equation}
and the average electric field is the convolution
\begin{equation}
\label{eq:Ebar}
\widebar{\V{E}}(\bm{z}) = \int_{\bm{x}} \bm{E}(\bm{x}) g\left(\norm{\V{x}-\V{z}}\right) \, d\bm{x}. 
\end{equation}
The force on each charge is given by
\begin{equation}
\label{eq:fqE}
\bm{F}_k = -\frac{\partial U}{\partial \V{z}_k} = q_k \widebar{\bm{E}}(\bm{z}_k). 
\end{equation}
Our goal will be to compute $U(\V{z})$ and $\V{F}(\V{z})$ to high accuracy in (log) linear time in the number of charges $N$. 


%

\section{Solver components \label{sec:bkgrnd}}
In this section, we build on prior work to introduce solvers for two simplified problems: the Poisson equation\ \eqref{eq:poiss} for smooth charge density $f$ and jump BCs\ \eqref{eq:cbc1}$-$\eqref{eq:cbc4}, and a free space solver for point charges. In Section\ \ref{sec:ourstuff}, we combine these two pieces to yield a solver for\ \eqref{eq:poiss}$-$\eqref{eq:cbc4} that remains efficient as $g_w \rightarrow 0$. 

We begin by introducing a solution method for \emph{smooth} doubly periodic problems, where a grid-based solver can efficiently resolve the charge density $f(\V{x})$, as would be the case for large $g_w$. Our method in this case is to use superposition to split the problem into two sub-problems: a Poisson equation over free space with uniform permittivity, and a harmonic correction to account for the jump conditions. This solver is inefficient for (near) point charges (small $g_w$) because the number of grid points necessary to resolve $f(\V{x})$ becomes too large. We address this in Section\ \ref{sec:ewnoBC} via Ewald splitting. 

\subsection{Smooth doubly periodic problems \label{sec:smoothdp}}
Let us suppose first that the distance between the charges is comparable to their width $g_w$. In this case, a grid-based method is an efficient way to solve\ \eqref{eq:poiss} with the BCs\ \eqref{eq:cbc1}$-$\eqref{eq:cbc4}. While this can be done with a single solve, it will aid us in Section\ \ref{sec:ourstuff} to split the potential into two pieces, 
\begin{equation}
    \phi=\phi^*+\cor{\phi}, 
\end{equation}
and solve for each piece separately. The first of these, $\phi^*$, is smooth and found by solving the Poisson equation\ \eqref{eq:poiss}, but with free space BCs and a \emph{uniform} permittivity $\epsilon$ (see\ \ref{sec:p1smooth}). The second piece, $\cor{\phi}$, is found by solving a Laplace equation for the correction potential that satisfies the BCs\ \eqref{eq:cbc1}$-$\eqref{eq:cbc4} (see\ \ref{sec:harmsolve}). This part of the solution is necessarily only piecewise smooth, since it is built to satisfy jump BCs. Because each of the pieces accounts for either the Gaussian or wall-bound charge densities (but not both), neither of them is necessarily electroneutral, and we therefore conclude in Section\ \ref{sec:combsmooth} by constructing a combined solution $\phi$ that gives a vanishing electric field as $z \rightarrow \pm \infty$.

\subsubsection{First problem: Poisson solve in free space \label{sec:p1smooth}}
For our first problem, we let $\phi^*$ be defined everywhere as the solution of a Poisson equation with uniform permittivity $\epsilon$ and no wall-bound charge densities, 
\begin{gather}
\label{eq:poiss2}
\epsilon \Delta \phi^*(\V{x}) = -\sum_{k=1}^M q_k g\left(\norm{\V{x}-\V{z}_k}\right). 
\end{gather}
We take free space boundary conditions in $z$, so that $\partial \phi/\partial z (z \rightarrow \pm \infty) \rightarrow 0$, and periodic boundary conditions in $x$ and $y$ on $[0,L_x] \times [0,L_y]$. For the moment we have no recourse to a grid-based solver, since the $z$ domain is unbounded. But since $f$ is a sum of charge densities concentrated on $(0,H)$, the Poisson equation\ \eqref{eq:poiss2} is actually a Laplace equation in the exterior, 
\begin{equation}
\label{eq:lapoutside}
    \epsilon \Delta \phi^* = 0 \text{ on $z < 0$ and $z > H$}. 
\end{equation}
This Laplace equation can be solved analytically by taking a Fourier transform in $x$ and $y$, which we denote with a hat. Using wave numbers $\kparv = (k_x,k_y)=(2\pi n/L_x, 2\pi m/L_y)$, with $m$ and $n$ integers and $\kpar = \norm{\kparv}$, we obtain
\begin{equation}
\frac{\partial ^2 \four{\phi}^*(\kparv,z)}{\partial z^2}-\kpar^2 \four{\phi}^*(\kparv,z) =0,
\end{equation}
which has the analytical solution (ruling out growth at infinity)
\begin{equation}
\label{eq:sollap}
    \four{\phi}^*(\kparv,z \leq 0) = C_2 e^{\kpar z} \;  \textrm{and} \; \four{\phi}^*(\kparv,z \geq H) = C_1 e^{-\kpar z},
\end{equation}
where $C_1$ and $C_2$ are unknown constants. 

The form of the solution\ \eqref{eq:sollap} implies the boundary conditions 
\begin{equation}
\label{eq:BCs1}
    \frac{\partial \four{\phi}^*(\kparv,z \rightarrow 0^-)}{\partial z}  -  \kpar \four{\phi}^*(\kparv,z \rightarrow 0^-)=0, \quad \frac{\partial \four{\phi}^*(\kparv,z \rightarrow H^+)}{\partial z}  +  \kpar \four{\phi}^*(\kparv,z\rightarrow H^+)=0.
\end{equation}
Since $\phi^*$ and $\partial \phi^*/\partial z$ are continuous everywhere, including at $z=0$ and $z=H$, the boundary conditions\ \eqref{eq:BCs1} must hold for the solution $\four{\phi}$ on $[0,H]$ as well. These equations are the \textit{Dirichlet to Neumann} map for a doubly-periodic domain \cite{DoublyPeriodic_Poisson, liang2020harmonic}, and they give us a two point boundary value problem (BVP) on $[0,H]$ for each Fourier mode, 
\begin{gather}
\label{eq:BVP}
\epsilon \left( \frac{\partial ^2 \four{\phi}^*(\kparv,z)}{\partial z^2}-\kpar^2 \four{\phi}^*(\kparv,z)\right) = -\four{f}(\kparv,z), \quad z \in [0,H],\\
\nonumber 
\frac{\partial \four{\phi}^*(\kparv,0)}{\partial z}  -  \kpar \four{\phi}^*(\kparv,0)=0, \quad \frac{\partial \four{\phi}^*(\kparv,H)}{\partial z}  +  \kpar \four{\phi}^*(\kparv,H)=0.
\end{gather}
Appendix \ref{sec:bvps} describes the integral formulation we use to numerically solve this BVP in the Chebyshev basis \cite{greengard1991spectral}. 

The solution for $\four{\phi}^*\left(\kparv=\V{0},z\right)$ is only well-defined when the integral of $f(\V{x})$ is zero. Because this is not in general the case, we will not be able to complete the formulation unless we consider both the charge densities $f(\V{x})$ and the surface charge densities $\sigma_{b/t}$ \emph{together}, which we do in Section\ \ref{sec:combsmooth}. For the moment, we consider the BVP for $\kparv=\V{0}$, 
\begin{equation}
\label{eq:k0eqn}
\epsilon \frac{\partial ^2 \four{\phi}^*(\V{0},z)}{\partial z^2}=-\four{f}(\V{0},z),
\end{equation}
and define an initial solution which is only correct up to a linear mode,  
\begin{equation}
\label{eq:phi0gen}
    \four{\phi}^*(\V{0},z) = -\frac{1}{\epsilon}\left(\int_{z'=0}^z \int_{s=0}^{z'}  \four{f}(0,s) \, ds \, dz'\right).
\end{equation}
The linear mode will be corrected in Section\ \ref{sec:combsmooth}.

\subsubsection{Second problem: harmonic correction \label{sec:harmsolve}}
We now move to the second piece of the potential in the slab, the correction $\cor{\phi}$, the purpose of which is to give the correct boundary conditions\ \eqref{eq:cbc1}$-$\eqref{eq:cbc4} for the total solution $\four{\phi}$. Because $\four{\phi}^*$ already included the charge densities in a Poisson solve, the only charge densities that remain are those on the walls. We therefore have the Laplace equation
 \begin{equation}
 \label{eq:harmeq}
 \Delta \cor{\phi}=0
 \end{equation}
 on $-\infty < z < 0$, $0 < z< H$, and $H < z<\infty$, all with $x$ and $y$ periodic. For generality, we assume that the boundary conditions are the jump conditions
 \begin{gather}
 \label{eq:harmbc1}
\cor{\phi}(x,y,z \rightarrow 0^+)-\cor{\phi}(x,y,z\rightarrow0^-) =-m^{(b)}_\phi(x,y),\\
\label{eq:harmbc2}
\epsilon \frac{\partial \cor{\phi}}{\partial z}(x,y,z\rightarrow0^+)-\epsilon_b \frac{\partial \cor{\phi}}{\partial z}(x,y,z\rightarrow0^-)=-m^{(b)}_E(x,y),\\
\label{eq:harmbc3}
\cor{\phi}(x,y,z\rightarrow H^-)-\cor{\phi}(x,y,z\rightarrow H^+) =-m^{(t)}_\phi(x,y),\\
 \label{eq:harmbc4}
\epsilon \frac{\partial \cor{\phi}}{\partial z}(x,y,z\rightarrow H^-)-\epsilon_t \frac{\partial \cor{\phi}}{\partial z}(x,y,z\rightarrow H^+)=-m^{(t)}_E(x,y). 
\end{gather}

Specifically, the values of $m_\phi$ and $m_E$ that ensure the overall BCs\ \eqref{eq:cbc1}$-$\eqref{eq:cbc4} are satisfied by $\phi=\phi^*+\cor{\phi}$ are given by
\begin{gather}
    m^{(b)}_\phi = \phi^*(x,y,z \rightarrow 0^+)-\phi^*(x,y,z \rightarrow 0^-)=0,\\
    m^{(b)}_E = \epsilon \frac{\partial \phi^*}{\partial z}(x,y,z \rightarrow 0^+)-\epsilon_b \frac{\partial \phi^*}{\partial z}(x,y,z \rightarrow 0^-)+\sigma_b(x,y)=(\epsilon-\epsilon_b)\frac{\partial \phi^*}{\partial z}(x,y,0)+\sigma_b(x,y),\\
    m^{(t)}_\phi = \phi^*(x,y,z \rightarrow H^-)-\phi^*(x,y,z \rightarrow H^+)=0,\\
    m^{(t)}_E = \epsilon \frac{\partial \phi^*}{\partial z}(x,y,z\rightarrow H^-)-\epsilon_t \frac{\partial \phi^*}{\partial z}(x,y,z \rightarrow H^+)-\sigma_t(x,y)=(\epsilon-\epsilon_t)\frac{\partial \phi^*}{\partial z}(x,y,H)-\sigma_t(x,y). 
\end{gather}
These particular expressions simplify because $\phi^*$ is continuously differentiable across $z=0$ and $z=H$, but in Section\ \ref{sec:ourstuff} we will have nonzero $m_\phi$.

The solution method for the harmonic problem\ \eqref{eq:harmeq}$-$\eqref{eq:harmbc4} is now straightforward. We introduce smooth harmonic functions $\cor{\phi}_i$, $\cor{\phi}_t$, and $\cor{\phi}_b$, whose domain is all of $\mathbb{R}^3$, and set
\begin{equation}
\label{eq:phicorgen}
    \cor{\phi} = \begin{cases} 
    \cor{\phi}_i & \text{if } z \in (0,H) \\ 
    \cor{\phi}_t & \text{if } z > H \\
    \cor{\phi}_b & \text{if } z < 0.
    \end{cases}
\end{equation}
After a Fourier transform in $x$ and $y$, we obtain the solution analytically as
\begin{gather}
\label{eq:harmfieldslab}
\cor{\four{\phi}}_i(\kparv,z )=\begin{cases} A_i(\kparv)e^{\kpar z}+B_i(\kparv)e^{-\kpar z} & \kpar > 0 \\[2 pt] A_i(\V{0})z+B_i(\V{0}) & \kpar =0, \end{cases}
\end{gather}
with solutions of the same form for $\four{\phi}_b(\kparv,z)$ and $\four{\phi}_t(\kparv,z)$, except that boundedness of $\four{\phi}$ implies that $A_t(\kpar \neq 0)=0$ and $B_b(\kpar \neq 0)=0$. When $\kpar \neq 0$, the four coefficients $A_{i/b}\left( \kparv \right)$, $B_{i/t}\left( \kparv \right)$ can be determined in a straightforward way from the four boundary conditions\ \eqref{eq:harmbc1}$-$\eqref{eq:harmbc4}. The resulting solution in the slab interior is
\begin{align}
\label{eq:inphisc}
\four{\phi}_i(\kparv,z )=&\frac{(\epsilon_t+1) e^{-\kpar z} (\four{m}^{(b)}_E(\kparv)-\epsilon_b \kpar
   \four{m}_\phi^{(b)}(\kparv))}{\kpar (\epsilon_b+1) (\epsilon_t+1) -(\epsilon_b
   \epsilon_t+\epsilon_b+\epsilon_t-1)e^{-2 H \kpar}}\\
   \nonumber
   -&\frac{(\epsilon_b+1) e^{\kpar (-H+z)} (\epsilon_t \kpar
   \four{m}_\phi^{(t)}(\kparv)+\four{m}_E^{(t)}(\kparv))}{\kpar (\epsilon_b+1) (\epsilon_t+1) -(\epsilon_b
   \epsilon_t+\epsilon_b+\epsilon_t-1)e^{-2 H \kpar}}\\
\nonumber
+&\frac{(\epsilon_b-1) e^{(-H-z) \kpar} (\epsilon_t \kpar
   \four{m}_\phi^{(t)}(\kparv)+\four{m}_E^{(t)}(\kparv))}{\kpar (\epsilon_b+1) (\epsilon_t+1) -(\epsilon_b
   \epsilon_t+\epsilon_b+\epsilon_t-1)e^{-2 H \kpar}}\\
   \nonumber
+&\frac{(\epsilon_t-1) \left(-e^{\kpar (z-2H)}\right) (\four{m}^{(b)}_E(\kparv)-\epsilon_b \kpar
   \four{m}_\phi^{(b)}(\kparv))}{\kpar (\epsilon_b+1) (\epsilon_t+1) -(\epsilon_b
   \epsilon_t+\epsilon_b+\epsilon_t-1)e^{-2 H \kpar}}
\end{align}
The solution for $\kpar=0$ is only well-defined when $\four{\phi}^*$ and $\cor{\four{\phi}}$ are added together, as we explain next.

\subsubsection{Electroneutrality and the $\kpar=0$ mode \label{sec:combsmooth}}
For the $\kparv=0$ mode, $\four{\phi}^*(\V{0},z)$ is defined as the double integral\ \eqref{eq:phi0gen}.
We also have linear modes for $\cor{\phi}_i=A_i(\V{0})z+B_i(\V{0})$ and likewise for $\cor{\phi}_b$ and $\cor{\phi}_t$. The total solution $\four{\phi}(\V{0},z) = \four{\phi}^*(\V{0},z)+\cor{\four{\phi}}(\V{0},z)$
must satisfy the slab BCs\ \eqref{eq:cbc1}$-$\eqref{eq:cbc4}. Substituting the representation\ \eqref{eq:harmfieldslab} for $\cor{\four{\phi}}$ into the slab BCs and using\ \eqref{eq:phi00}, we obtain equations for the unknown coefficients in terms of the double integral $\four{\phi}^*(\V{0},z)$, 
\begin{gather}
\label{eq:p0tot}
    \four{\phi}^*(\V{0},0)+B_i(\V{0}) = B_b(\V{0})=0,\\
    \label{eq:p2tot}
    A_i(\V{0})H+B_i(\V{0})+\four{\phi}^*(\V{0},H) = A_t(\V{0})H + B_t(\V{0}),\\
    \label{eq:E1tot}
    \epsilon \left(\frac{\partial \four{\phi}^*}{\partial z} (\V{0},0)+A_i(\V{0})\right) - \epsilon_b\left(\frac{\partial \four{\phi}^*}{\partial z} (\V{0},0)+ A_b(\V{0})\right)+\four{\sigma}_b(\V{0})=0,\\
    \label{eq:E4tot}
    \epsilon \left(\frac{\partial \four{\phi}^*}{\partial z} (\V{0},H)+A_i(\V{0})\right) - \epsilon_t \left(\frac{\partial \four{\phi}^*}{\partial z} (\V{0},H)+A_t(\V{0})\right)-\four{\sigma}_t(\V{0})=0.
\end{gather}
There are two more equations relating to the decay of the electric field at $z \rightarrow \pm \infty$. On $z < 0$ and $z > H$, the electric field must decay to zero, which implies that the growth from $\cor{\phi}_b$ and $\cor{\phi}_t$ must cancel that due to $\four{\phi}^*$ on $z \leq 0$ and $z \geq H$, 
\begin{equation}
\label{eq:Atb}
    A_b(\V{0}) = -\frac{\partial \four{\phi}^*}{\partial z} (\V{0},z \leq 0) \qquad  A_t(\V{0}) = -\frac{\partial \four{\phi}^*}{\partial z} (\V{0}, z \geq H).
\end{equation}
Substituting this into\ \eqref{eq:E1tot}$-$\eqref{eq:E4tot}, we obtain two solutions for $A_i(\V{0})$, which are equivalent for electroneutral slabs,  
\begin{align}
\label{eq:Ai1}
   A_i (\V{0}) & = -\frac{\four{\sigma}_b(\V{0})}{\epsilon}-\frac{\partial \four{\phi}^*}{\partial z} (\V{0},z = 0)\\
   \label{eq:Ai2}
   & = \frac{\four{\sigma}_t(\V{0})}{\epsilon}-\frac{\partial \four{\phi}^*}{\partial z} (\V{0},z = H).
\end{align}
The equivalence can be seen by applying the divergence theorem to the Poisson equation\ \eqref{eq:poiss2}, and then replacing the surface integrals with zero-Fourier modes to obtain
\begin{equation}
    \frac{\partial \four{\phi}^*}{\partial z} (\V{0},z = H)- \frac{\partial \four{\phi}^*}{\partial z}(\V{0},z = 0) = -\frac{1}{\epsilon} \sum_{k=1}^N q_k.
\end{equation}
Substituting this into the electroneutrality condition\ \eqref{eq:elecneu} yields\ \eqref{eq:Ai1}$=$\eqref{eq:Ai2}. The values of $B_i(\V{0})$, $B_b(\V{0})$, and $B_t(\V{0})$ are determined from\ \eqref{eq:p0tot} and\ \eqref{eq:p2tot}.

\subsection{Ewald splitting \label{sec:ewnoBC}}
Accurately computing the potential using the method of Section\ \ref{sec:smoothdp} for a collection of Gaussian charges requires a grid spacing $h \sim g_w$. If $g_w$ is comparable to the distance between charges, then this method is sufficient to solve the problem. But as $g_w \rightarrow 0$ (e.g., the case of point charges), the method becomes woefully inefficient as the grid size is dictated by the size of the charges and not the spacing between them. 

In this section, we address this problem by Ewald splitting \cite{SpectralEwald_Electrostatics, SpectralEwaldElectrostatics_Review}. The idea of Ewald splitting or Ewald summation is to smear the charges to the extent that they can be resolved on a grid of reasonable size, and then correct for this smearing through an analytical near-field correction. Here we present an Ewald splitting algorithm \emph{without} dielectric jumps or periodicity \cite{SpectralEwald_Electrostatics, SpectralEwaldElectrostatics_Review}. We then specialize to the doubly-periodic case with dielectric interfaces in Section\ \ref{sec:ourstuff}. 

As discussed in the introduction, previous presentations of SE methods \cite{SpectralEwald_NFFT,SpectralEwald_Electrostatics,DoublyPeriodic_Electrostatics} view Ewald splitting as an algebraic construction for the fast calculation of electrostatic sums. The sums are split into two pieces, with the near-field sums computed by summing over neighboring points, and the far-field sums computed using non-uniform FFT methods \cite{NUFFT}. In the approach of \cite{SpectralEwald_Electrostatics}, this leads to two separate Gaussian kernels: one for the algebraic splitting, and one in the NUFFT method for spreading and interpolating to/from the FFT grid. This is advantageous in triply periodic domains, since the NUFFT-based approach allows for the use of any NUFFT method \cite{SpectralEwald_NFFT,NUFFT}, including recent methods based on non-Gaussian kernels \cite{FINUFFT_Barnett}. \rev{Extending this NUFFT-based approach to systems with mixed periodicity requires either expensive oversampling \cite{DoublyPeriodic_Electrostatics} or (non-analytic) function extension in the aperiodic directions \cite{SpectralEwald_NFFT}.} 

Since we instead view Ewald splitting as a smearing of charges to allow for a grid-based method, the only Gaussian kernel we use is the one in the smearing step, and this same kernel is used both in the Ewald splitting \emph{and} in the communication between particles and the grid-based solver (spreading and interpolation). For triply periodic domains, this gives an approach similar to that of \cite{SpectralEwald_Electrostatics}.\footnote{Specifically, our approach for triply periodic domains is equivalent to the SE method of \cite{SpectralEwald_Electrostatics} with the parameter $\eta$ fixed to $\eta=1$.} \rev{However, the smearing approach we take here straightforwardly generalizes to systems with mixed periodicity because our far-field Ewald sum only requires a method to solve a smooth Poisson equation in the domain of interest. That said, the DtN approach that we use here to solve the Poisson equation in 2D periodic systems does not efficiently extend to 1D periodic systems (e.g., a square channel) or 0D periodic systems (e.g., a cubic chamber) because the DtN map is nontrivial to compute (e.g., it is in general nonlocal).}

Beginning with the Poisson equation\ \eqref{eq:poiss}, posed in all of $\mathbb{R}^3$, we add and subtract a convolution (denoted by $\star$) with a ``screening'' or ``splitting'' function $\gamma$. We will do the splitting symmetrically by writing 
\begin{gather}
\label{eq:lapsplit}
\phi = \underbrace{\gamma^{1/2} \star \psi}_{\far{\phi}}+ \underbrace{\left(\phi - \gamma^{1/2} \star \psi \right)}_{\near{\phi}} , \\
\label{eq:fmod}
\text{where } \epsilon \Delta \psi(\V{x}) = -\left(f \star \gamma^{1/2}\right)(\V{x}). 
\end{gather}
The first term $\far{\phi}$ in\ \eqref{eq:lapsplit} is the ``far field'' and is the potential due to smeared versions of the original charges. We smear the Gaussian charge densities\ \eqref{eq:gcharge} using a radially symmetric Gaussian kernel with standard deviation $1/(2\xi)$ \cite{SpectralEwald_Electrostatics},
\begin{equation}
\label{eq:gamreal}
\gamma^{1/2}(r;\xi) = \frac{8\xi^3 }{(2\pi) ^{3/2}}e^{-2r^2 \xi^2},
\end{equation}
where the ``splitting parameter'' $\xi$ is an arbitrary positive real number that is chosen to optimize speed (see Section \ref{sec:performance}). 

Although the far field potential defined in\ \eqref{eq:lapsplit} is the convolution of $\gamma^{1/2}$ with $\psi$, we recall that we seek the potential averaged over the Gaussian cloud $g$, so that our real goal is to obtain
\begin{align}
\label{eq:compGfar}
\far{\widebar{\phi}}=g\star\far{\phi}=\left(g\star \gamma^{1/2}\right) \star \psi. 
\end{align}
The kernel $\widebar{S}:=g\star \gamma^{1/2}$ is a convolution of two Gaussians, one with standard deviation $g_w$ and the other with standard deviation $1/(2\xi)$,
\begin{gather}
\label{eq:spreadkernelEw}
\widebar{S}(r;g_w,\xi)  = g\star\gamma^{1/2}=\frac{1}{\sqrt{8\pi^3 g_t^6}} \exp{\left(-\frac{r^2}{2g_t^2}\right)},\\
\label{eq:gt} \text{where} \quad g_t = \sqrt{\frac{1}{4\xi^2}+g_w^2}  \, ,
\end{gather}
that is, $\widebar{S}$ is a Gaussian with standard deviation $g_t$. We can now write the average far field potential $\far{\phi}=\widebar{S} \star \psi$ in the symmetric form
\begin{equation}
    \far{\phi}(\V{x}) = \widebar{S} \star \Delta^{-1} \left(\sum_{k=1}^N q_k \widebar{S}\left(\norm{\V{x}-\V{z}_k}\right)\right). 
\end{equation}
Here $\psi(\V{x})$ can be computed on a grid with spacing $h \sim g_t$ using any grid-based solver, for example an FFT-based solver for triply periodic domains. 

The remaining term $\near{\phi}$ in\ \eqref{eq:lapsplit} is the ``near field.'' Since the Laplacian commutes with the convolution, we use\ \eqref{eq:fmod} to write a Poisson equation for the near field
\begin{equation}
\label{eq:nearprob}
\epsilon \Delta \near{\phi} = -f + \gamma^{1/2} \star \left(f \star \gamma^{1/2}\right)= -f \star(1-\gamma).
\end{equation} 
The ``near field'' charge density $f \star (1-\gamma)$ for a single charge is sharply peaked, has negative tails, and integrates to zero in three dimensions. Thus the field it creates decays rapidly (exponentially) in real space, and can be computed analytically using Fourier integrals. For this, we recall that the charge density $f$ for a single charge is proportional to the Gaussian $g(r)$ defined in\ \eqref{eq:gcharge}, with Fourier transform 
\begin{equation}
\label{eq:gaussFour}
\four{g}(k) = \exp{\left(-\frac{1}{2}g_w^2 k^2\right)}.
\end{equation}
We obtain a near field interaction kernel $\near{G}$ by setting $f=g$ in\ \eqref{eq:nearprob} and solving the resulting algebraic equation in Fourier space, 
\begin{gather}
\label{eq:Gnearfour}
\near{\four{G}}(k;g_w,\xi) = \frac{\left(1-\four{\gamma}\right)\four{g}(k)}{\epsilon k^2}, \quad \text{where}\\
\four{\gamma}(k;\xi) = \exp{\left(-k^2/4\xi^2\right)}.
\end{gather}
The near field interaction kernel in real space can now be obtained by a radially-symmetric inverse Fourier transform of\ \eqref{eq:Gnearfour},  
\begin{align}
\nonumber
\near{G}(r;g_w,\xi) & = \frac{1}{2\pi^2r} \int_0^\infty  \frac{\left(1-\four{\gamma}(k)\right)\four{g}(k)}{\epsilon k^2} k \sin{(kr)} \, dk \\[4 pt]
\label{eq:Gnear}
& = \frac{1}{4\pi \epsilon r} \left(\textrm{erf}{\left(\frac{r}{\sqrt{2}g_w}\right)-\textrm{erf}{\left(\frac{r}{\sqrt{2g_w^2+\xi^{-2}}}\right)}}\right). 
\end{align}
 In the limit as $\xi \rightarrow 0$ and $g_w \rightarrow 0$, we recover the equation for the potential of a point charge in free space, $1/(4\pi\epsilon r)$, while in the limit $g_w \rightarrow 0$ with $\xi$ fixed, we obtain the near field interaction kernel for point charges \cite{SpectralEwaldElectrostatics_Review}. Because the near field interaction kernel\ \eqref{eq:Gnear} decays exponentially, it can be truncated to zero for $r > \rnf$, where $\rnf \sim g_t$ (see Section\ \ref{sec:nearconstraint}). 
 
 The near field kernel\ \eqref{eq:Gnear} will give the pointwise near field potential at any point in space via
\begin{equation}
\label{eq:Gnearsum}
    \near{\phi}(\V{x}) = \sum_{k=1}^N q_k \near{G}(\norm{\V{x}-\V{z}_k};g_w,\xi)
\end{equation}
 The \emph{average} near field potential $\near{\widebar{\phi}}=g \star \near{\phi}$ is given by 
\begin{equation}
\label{eq:avgGn}
    \near{\widebar{\phi}}(\V{x}) =\sum_{k=1}^N q_k\near{\widebar{G}}(\norm{\V{x}-\V{z}_k};g_w,\xi)
\end{equation}
where the \emph{average} near field interaction kernel, $\near{\widebar{G}}=g\star\near{G}$. This convolution can be computed analytically, since $\near{\four{\widebar{G}}}=\four{g} \near{\four{G}}$, to give
\begin{gather}
\label{eq:GnearK}
\near{\widebar{G}}(r;g_w,\xi) = \frac{1}{4\pi \epsilon r} \left(\textrm{erf}{\left(\frac{r}{2g_w}\right)-\textrm{erf}{\left(\frac{r}{\sqrt{4g_w^2+{\xi^{-2}}}}\right)}}\right) = \near{G}(r;g_w\sqrt{2},\xi). 
\end{gather}

In a similar way to the potential field, we define the pointwise near field electric field by
\begin{equation}
\label{eq:Enear}
\near{\V{E}}(\V{x};g_w,\xi) = - \sum_{k=1}^N q_k \frac{\partial \near{G}(\norm{\V{r}_k};g_w,\xi)}{\partial r} \widehat{\V{r}}_k, 
\end{equation}
where $\V{r}_k=\V{x}-\V{z}_k$ and $\widehat{\V{r}}_k=\V{r}_k/\norm{\V{r}_k}$. The averaged near field electric field is given by convolving\ \eqref{eq:Enear} with $g$, or, equivalently, by differentiating the average near field Green's function\ \eqref{eq:GnearK},
\begin{equation}
\label{eq:EnearK}
\near{\widebar{\V{E}}}(\V{x};g_w,\xi)
= - \sum_{k=1}^N q_k \frac{\partial \near{\widebar{G}}(\norm{\V{r}_k};g_w,\xi)}{\partial r} \widehat{\V{r}}_k.  
\end{equation}
To avoid cancellation of digits, for small $r$ values ($r < 10^{-3}g_w$ in double precision) we use the Taylor series for $\near{\widebar{G}}$ in\ \eqref{eq:EnearK}. 

\section{Ewald splitting for slabs \label{sec:ourstuff}}
In this section we develop a solver for point-like charges in a slab geometry. The first step is to define the solution using \emph{images} rather than boundary conditions at the slab walls. For a single dielectric boundary, reflecting the charge locations across the boundary and giving them modified strengths yields a set of image charges. The potential from these images plus the original charges is then what we seek. We show in Section\ \ref{sec:imgs} that this method of images can be applied to the near field and far field problems separately, with the result that both fields satisfy the boundary conditions at the dielectric interface. Thus for a single dielectric jump, we could make a set of images and solve the problem with uniform permittivity using Ewald splitting as described in Section\ \ref{sec:ewnoBC}.

The situation for multiple dielectric jumps is more complex because there are infinitely many images. The novelty in our algorithm is its ability to handle the infinitely many images in an efficient way. Beginning with the system of infinitely many images, we follow the Ewald splitting method of Section\ \ref{sec:ewnoBC} to form a near field and far field problem with the infinitely many image charges. In the near field problem, the interaction kernel\ \eqref{eq:GnearK} decays exponentially in real space, and therefore in Section\ \ref{sec:ournear} we truncate it so that the charges inside the slab only interact with their nearest images on either side of the slab. 

In the far field problem, the smearing of the charges and images allows for the use of a grid-based solver, and in Section\ \ref{sec:farfield} we evoke the components developed in Section\ \ref{sec:smoothdp} to split the far field problem\ \eqref{eq:fmod} into two solves. In the first solve, described in Section\ \ref{sec:intpots}, we construct an initial solution $\psi^*$ for the potential inside and outside the slab. We include in this initial solve a requisite number of images so that the remaining field $\cor{\psi}=\psi-\psi^*$ is harmonic and can be determined analytically in $(\kparv,z)$ space as described in Section\ \ref{sec:farharm}. As in Section\ \ref{sec:smoothdp}, the $\kpar=0$ mode is only well-defined for the total potential $\psi=\psi^*+\cor{\psi}$, and so we consider it at the end in Section\ \ref{sec:farcor}. 

\subsection{Image construction \label{sec:imgs}}
We begin by reviewing the classical image construction for point charges near a single dielectric interface at $z=0$, where for $z > 0$ the permittivity is $\epsilon$ and for $z < 0$ the permittivity is $\epsilon_b$. Let a charge of strength $q$ be centered at $(x,y,z>0)$. Then the potential and electric field for $z > 0$ are as if the medium had uniform permittivity $\epsilon$ and contained the original charge along with an image charge of strength
\begin{equation}
\label{eq:ims}
q^{*} = -q \frac{\epsilon_b-\epsilon}{\epsilon_b+\epsilon},
\end{equation}
positioned at $(x,y,-z)$. For $z < 0$, the potential and electric field are again as if the medium had uniform permittivity $\epsilon$ and contained a single charge positioned at $(x,y,z)$, with strength
\begin{equation}
\label{eq:imb}
q^{**} = q \frac{2\epsilon}{\epsilon_b+\epsilon}.
\end{equation}
This image construction ensures that the pointwise potential $\phi$ and the normal component of the electric displacement $\epsilon \V{E}$ are continuous at the wall, i.e., that the boundary conditions\ \eqref{eq:cbc1} and\ \eqref{eq:cbc3} are satisfied with $\sigma_b=0$. 

The same image construction can be used when the charges are spherically-symmetric clouds instead of points. To show this, let $G(r)$ be the function that gives the pointwise potential a distance $r$ from the center of the charge. We assume that the charge is centered at height $\Delta z$ above the interface.  Let $(x,y,0)$ be any point on the $xy$ plane, and let $r$ be the distance from the center of the charge to that point. Then, using the image construction, the potential at $z=0$ from $z < 0$ is given by (using\ \eqref{eq:imb})
\begin{equation}
\phi(z \rightarrow 0^-) = \frac{2\epsilon}{\epsilon_b+\epsilon}G(r),
\end{equation}
while using\ \eqref{eq:ims} for the potential on $z> 0$ we get
\begin{equation}
\phi(z \rightarrow 0^+)  = 
G(r) - \frac{\epsilon_b-\epsilon}{\epsilon_b+\epsilon} G(r) = \phi(z \rightarrow 0^-).
\end{equation}
The electric displacement calculation is similar. From\ \eqref{eq:imb}, we have 
\begin{equation}
\epsilon_b \frac{ \partial \phi}{\partial z}(z \rightarrow 0^-) = \frac{2\epsilon \epsilon_b}{\epsilon_b+\epsilon}\frac{\partial G}{\partial r}(r)\frac{\Delta z}{r},
\end{equation}
while using\ \eqref{eq:ims} we obtain
\begin{gather}
\epsilon \frac{ \partial \phi}{\partial z}(z \rightarrow 0^+) = \epsilon  \frac{\partial G}{\partial r}(r)\left(\frac{\Delta z}{r}- \left(\frac{\epsilon_b-\epsilon}{\epsilon_b+\epsilon}\right) \frac{(-\Delta z)}{r}\right) = \epsilon_b \frac{ \partial \phi}{\partial z}(z \rightarrow 0^-),
\end{gather}
which confirms that the potential and electric displacement are continuous across the interface. Since the form of the kernel $G(r)$ does not matter, this implies that the classical image construction applies to any isotropic charge clouds (in particular, Gaussian clouds), even ones that overlap the interface.

\subsection{Near field in the slab geometry \label{sec:ournear}}
For the doubly periodic slab geometry, the rapid decay of the near field kernel can be used to simplify the number of required images in the near field problem. Let us consider first the case of a single wall to simplify the argument. Recall that $\rnf$ is the distance at which the near field kernel is truncated. Then in order for the near field pointwise potential $\near{\phi}$ to satisfy the boundary conditions at the wall, we must include all images that are centered $\rnf$ or closer to the wall. 

\begin{figure}
\centering     
\includegraphics[width=0.7\textwidth]{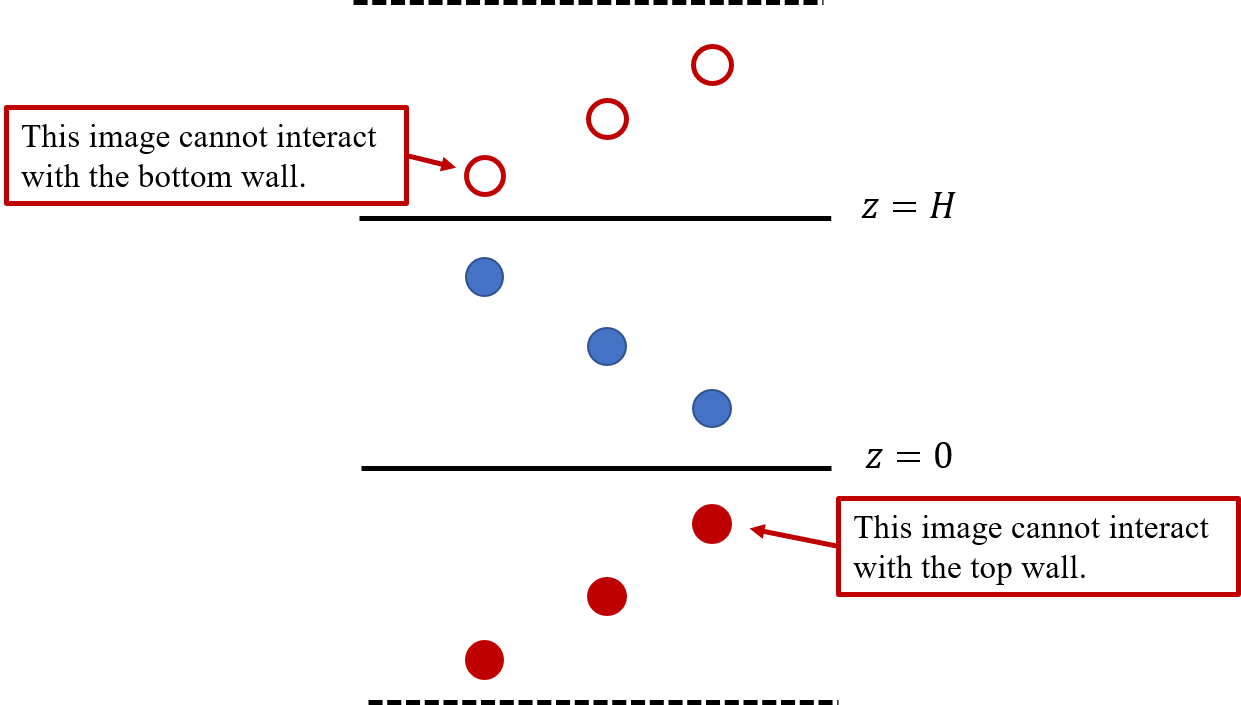}
\caption{\label{fig:nearfieldimg} Images of the original charges in the slab (blue circles) in the near field problem. If we want the jump BCs\ \eqref{eq:cbc1}$-$\eqref{eq:cbc4} to be satisfied with $\sigma_t=\sigma_b=0$ and only one set of images above/below the slab, then the near field created by the images below $z=0$ (filled red circles) cannot extend above $z \geq H$, and similarly for the images above $z > H$ (empty red circles) and the bottom wall.}
\end{figure}

The case of the slab is similar. In order for the near field pointwise potential $\near{\phi}$ to satisfy boundary conditions at a wall, we only need to include the images (with respect to that wall) that are $\rnf$ or closer to it. Suppose that we want to limit the near field computation to the first set of images above/below each wall (we analyze the implications of this choice in Section\ \ref{sec:nearconstraint}). Then, as shown in Fig.\ \ref{fig:nearfieldimg}, the key requirement is that images below the bottom wall cannot affect the top wall. Now let $h$ be the minimum distance between the bottom wall of the slab and the center of any charge. Then the closest an image can be to the bottom wall is $h$, and we have the requirement that $\rnf < H+h$ for only the first set of images below $z < 0$ to be included in the near field sum. Similar considerations apply to the images above $z > H$ and the bottom wall. We also ensure that $\rnf < \min{\left(L_x/2,L_y/2\right)}$ so that we can apply the nearest image convention in the near field. 

Including only the first set of images is equivalent to treating the problem as two single walls, rather than as a full slab geometry. The slab geometry only enters the near field when we include more than the first image. By considering a picture similar to Fig.\ \ref{fig:nearfieldimg}, it can be shown that restricting the number of near field images to $n_\text{img}$ gives the requirement $\rnf < n_\text{img} H +h$. However, we will show in Section\ \ref{sec:xi} that far-field constraints on the Ewald parameter $\xi$ give an algorithm for which $\rnf < H +h$ automatically, and so we will set $n_\text{img}=1$. 

\subsection{Far field in the slab geometry \label{sec:farfield}}
We are now ready to tackle the far field problem. For the dielectric slab, we need to solve the far field Poisson equation\ \eqref{eq:fmod} with periodicity in the $x$ and $y$ directions and slab boundary conditions\ \eqref{eq:cbc1}$-$\eqref{eq:cbc4}. To do this, we define $\psi$ piecewise by 
\begin{equation}
\label{eq:psi}
    \psi(x,y,z) = \begin{cases}
    \psi_i(x,y,z) & 0 \leq z \leq H \\
    \psi_b(x,y,z) & z < 0 \\
    \psi_t(x,y,z) & z > H 
    \end{cases}, 
\end{equation}
where $\psi_i$, $\psi_b$, and $\psi_t$ are smooth fields defined for all $z \in \mathbb{R}$.\footnote{\rev{Note that $\psi_{i/t/b}$ are analytic (and thus $C^\infty$), but $\psi$ itself is not. Since we numerically only solve for $\psi_{i/t/b}$ and not for $\psi$ itself, our approach still gives spectral accuracy despite the lack of smoothness at the dielectric boundaries.}} Specifically, the potentials $\psi_{b}$ and $\psi_t$ are due to the far-field charges and their (infinitely many) images above $z =H$ and below $z=0$, respectively. Using $\farf:=f*\gamma^{1/2}$ to denote the charge density from the (smeared) far-field charges, $\img{\farf}_{z>H}$ to denote the density from the smeared images centered above the top wall, and $\img{\farf}_{z<0}$ to denote the density from the smeared images centered below the bottom wall, the potentials $\psi_b$ and $\psi_t$ satisfy the Poisson equations
\begin{gather}
\label{eq:psib}
\epsilon \Delta \psi_b = -\frac{2\epsilon}{\epsilon_b+\epsilon}\left(\farf+\img{\farf}_{z>H} \right),\\
\label{eq:psit}
\epsilon \Delta \psi_t = -\frac{2\epsilon}{\epsilon_t+\epsilon}\left(\farf+\img{\farf}_{z<0} \right),
\end{gather}
with periodic BCs in the $x$ and $y$ directions and free space BCs in the $z$ direction (i.e., no dielectric boundaries and decay of the potential as $z \rightarrow \pm \infty$). 

The interior potential $\psi_i$ satisfies the Poisson equation
\begin{equation}
\label{eq:farinf}
\epsilon \Delta \psi_i = -\left(\farf + \img{\farf} \right),
\end{equation}
with periodicity in the $x$ and $y$ directions. Here $\img{\farf}=\img{\farf}_{z>H}+\img{\farf}_{z<0}$ is the charge density due to the (infinite number of) images above and below the slab. Note that the strength of these images decays to zero as $z \rightarrow \infty$, so the Poisson equations\ \eqref{eq:psib}$-$\eqref{eq:farinf} are well-posed. If we use free space boundary conditions in the $z$ direction for\ \eqref{eq:farinf}, then we know that $\psi$ as defined in\ \eqref{eq:psi} will satisfy the jump BCs\ \eqref{eq:cbc1}$-$\eqref{eq:cbc4} with $\sigma_b=\sigma_t=0$ (i.e., with no wall-bound charge). 

Assuming that the densities $\sigma_{b/t}$ are smooth enough to be resolved by the grid used for the far field solve, it makes sense to include the potential from the wall-bound densities, which we denote by $\cor{\phi}$, as part of the far field potential $\far{\phi}$. On the slab interior, we define the far field potential as
\begin{equation}
    \far{\phi} := \gamma^{1/2} * \psi_i + \cor{\phi}_i, 
\end{equation}
where $\cor{\phi}_i$ is the harmonic potential defined in\ \eqref{eq:phicorgen} and\ \eqref{eq:harmfieldslab} with $m_\phi=0$ and $m_E= \pm \sigma_{b/t}$. 

Because $\cor{\phi}_i$ is harmonic, by the generalized mean value theorem it is indistinguishable from its convolution with the radially-symmetric kernel $\gamma^{1/2}$. This means that the far field potential can be rewritten as
\begin{equation}
    \far{\phi} = \gamma^{1/2} * \left(\psi_i + \cor{\phi}_i\right), \quad \text{and} \quad \far{\widebar{\phi}} = \widebar{S} * \left(\psi_i + \cor{\phi}_i\right).
\end{equation}
For convenience of notation, we redefine $\psi_i$ to refer to the combined potential $\psi_i+\cor{\psi}_i$ from the smeared charges and wall-bound densities. Using this definition, we obtain boundary conditions for\ \eqref{eq:farinf},
\begin{gather}
\label{eq:psiBC1}
\psi_i(x,y,z=0)-\psi_b(x,y,z=0)=0,\\
\label{eq:psiBC2}
\epsilon \frac{\partial \psi_i}{\partial z}(x,y,z=0)-\epsilon_b \frac{\partial \psi_b}{\partial z}(x,y,z=0)=-\sigma_b(x,y),\\
\label{eq:psibc3}
\psi_i(x,y,z=H)-\psi_t(x,y,z=H)=0,\\
\label{eq:psiBC4}
\epsilon \frac{\partial \psi_i}{\partial z}(x,y,z=H)-\epsilon_t \frac{\partial \psi_t}{\partial z}(x,y,z=H)=\sigma_t(x,y).
\end{gather}

\subsubsection{Potential and charge splitting}
\begin{figure}
\centering
\includegraphics[width=0.7\textwidth]{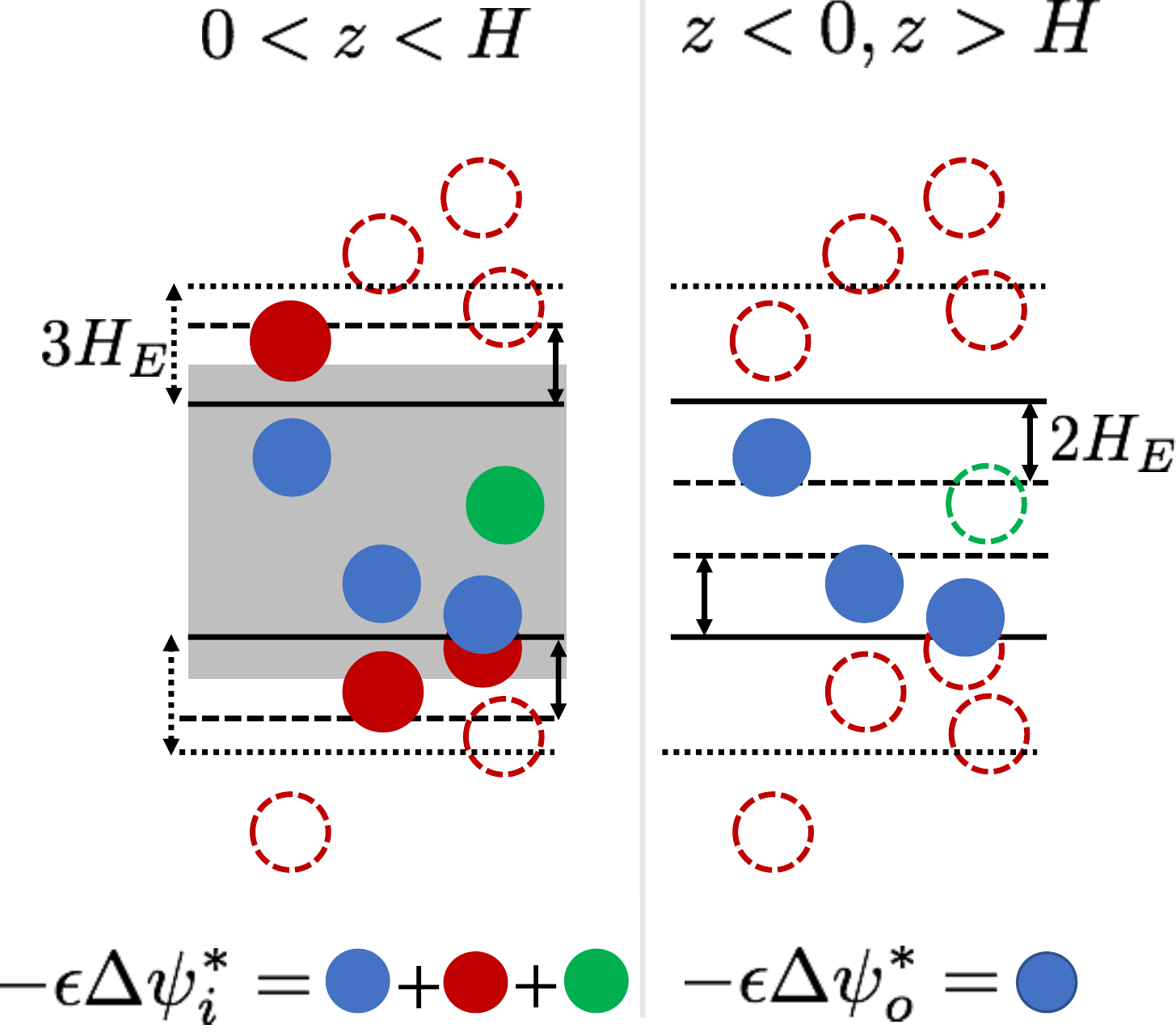}
\caption{\label{fig:ffsolve} Far field algorithm. We split the charges and images into four groups: blue charges = $\olap{C}$, green charges = $\nolap{C}$, filled red images = $C^\text{(img)}_{<2H_E}$, dashed red circles = $C^\text{(img)}_{>2H_E}$; the group $C^\text{(img)}_{>2H_E}$ contains infinitely many images, not all of which are shown. Circles denote the support of the smeared far-field charges, solid black lines denote the dielectric interfaces, solid arrows denote a region of size $2H_E$ (the diameter of the far-field charge support), and dotted arrows denote a region of size $3H_E$ (the amount we must extend the domain to enclose the solid circles).  Beginning in the left column, we first obtain an initial guess for the solution inside the slab, $\psi_i^*$, using the charges (solid blue and solid green) and images (solid red) whose support overlaps the gray region from which we interpolate, $ z \in [-H_E,H+H_E]$. In the right column, we obtain an initial guess for the solution outside the slab, $\psi_o^*$, using the charges (solid blue) that are centered $2H_E$ or closer to either wall. Both Poisson problems (for $\psi_i^*$ and $\psi_o^*$) are solved on a doubly periodic domain with $z \in [-3H_E,H+3H_E]$ using the method of Section\ \ref{sec:p1smooth}. The rest of the charges and images (dashed circles) can be coarse-grained into a harmonic correction solve for $\cor{\psi}_i$ and $\cor{\psi}_{t/b}$. }
\end{figure}

While the Poisson equations\ \eqref{eq:psib}$-$\eqref{eq:farinf} are periodic in $x$ and $y$ and unbounded in $z$, we cannot directly apply the techniques of Section\ \ref{sec:smoothdp} to solve them because the image charge densities in\ \eqref{eq:psib}$-$\eqref{eq:farinf} are unbounded in the $z$ direction. To work around this, we split each potential 
\begin{equation}
\label{eq:splitpsi}
    \psi_{i/t/b}=\psi_{i/t/b}^* + \cor{\psi}_{i/t/b}
\end{equation}
into the potential due to a finite number of images close to the slab (denoted with $*$ and obtained by the techniques of Section\ \ref{sec:p1smooth}), and the harmonic potential due to the rest of the (infinitely many) images (denoted with $(c)$ and obtained using the method of Section\ \ref{sec:harmsolve}). A conceptual picture of our far field solver is depicted in Fig.\ \ref{fig:ffsolve}.  

To more precisely quantify the domain where we need to evaluate $\psi_i$, we define the truncation radius of the smeared (far-field) Gaussian charges as $H_E$, where typically $H_E>4g_t$. If the charge centers are contained on $z \in [0,H]$, then we need to compute $\psi_i$ for $z \in [-H_E,H+H_E]$ so that we can average it with the kernel $\widebar{S}$ to obtain $\widebar{\psi}_i = \widebar{S} \star \psi_i$. Thus our goal is to obtain $\psi_i$ for at least $z \in [-H_E,H+H_E]$ (gray region in Fig.\ \ref{fig:ffsolve}).

We first split the smeared charge density $\farf$ into a piece $\nolap{\farf}$ coming from charges centered at least $2H_E$ from the slab boundaries (i.e., charges with centers in $z \in [2H_E,H-2H_E]$, denoted by $\nolap{C}$ and shown with green circles in Fig.\ \ref{fig:ffsolve}), and a piece $\olap{\farf}$ coming from charges centered a distance $2 H_E$ or less from the slab boundaries (denoted $\olap{C}$ and shown with blue circles in Fig.\ \ref{fig:ffsolve}), 
\begin{equation}
\label{eq:fsp}
\farf = \nolap{\farf} + \olap{\farf}.
\end{equation}
The splitting on the set of images $\img{\farf}=\nolap{\farf}^\text{(img)}+\olap{\farf}^\text{img}$ is exactly the same, with $\olap{\farf}$ representing the first set of images of $\olap{C}$, denoted by $\olap{C}^\text{(img)}$ and shown with solid red circles in Fig.\ \ref{fig:ffsolve}, and $\nolap{\farf}$ representing the rest of the images (including higher-order images of $\olap{C}$), denoted by $\nolap{C}^\text{(img)}$ and shown with dashed red circles in Fig.\ \ref{fig:ffsolve}.

Using this splitting, the Poisson equation\ \eqref{eq:farinf} for the interior solution can be written as
\begin{gather}
\label{eq:phii}
\epsilon \Delta \psi_i = -\left( \olap{\farf}+\nolap{\farf} +\olap{\farf}^\text{(img)}+\nolap{\farf}^\text{(img)}\right),
\end{gather}
 with the boundary conditions\ \eqref{eq:psiBC1}$-$\eqref{eq:psiBC4}. In Fig.\ \ref{fig:ffsolve}, we show (some of) the infinite number of images that make up the density $\nolap{\farf}^\text{(img)}$ as dashed red circles. We see that all of these images are centered at least $2H_E$ from the slab (i.e., in the region $z \in (\infty,-2H_E) \cup (H+2H_E,\infty)$ that is outside of the dashed lines in Fig.\ \ref{fig:ffsolve}), and are therefore compactly supported outside of the gray region $z \in [-H_E,H+H_E]$, on which we need to compute $\psi_i$. 

\subsubsection{Intermediate potentials \label{sec:intpots}}
We define the intermediate potential $\psi_i^*$ as the potential from all charges inside the slab, together with the (finite number of) images $\olap{C}^\text{(img)}$ whose support penetrates the region $z \in [-H_E,H+H_E]$. This is the potential due to the solid circles in the left column of Fig.\ \ref{fig:ffsolve} and is given by
\begin{gather}
\label{eq:phiistar}
\epsilon \Delta \psi_i^* = -\left(\nolap{\farf} + \olap{\farf}+\olap{\farf}^\text{(img)} \right),
\end{gather}
with periodicity in the $x$ and $y$ directions and free space BCs in the $z$ direction.  This Poisson equation can be solved using the method of Section\ \ref{sec:p1smooth}, except that the domain in the $z$ direction has to be modified. Because we include images centered $2H_E$ or closer to either wall, and because all smeared far-field charges have a support of $H_E$ in every direction, the right hand side (r.h.s.)\ of\ \eqref{eq:phiistar} is supported on $z \in [-3H_E,H+3H_E]$. Since all images with support overlapping the region $z \in [-H_E,H+H_E]$ are included in\ \eqref{eq:phiistar}, the difference $\cor{\psi}_i=\psi_i-\psi_i^*$ must be harmonic for $z \in [-H_E,H+H_E]$. 

Similarly, we construct solutions for $\psi_b$ and $\psi_t$ by using the splitting\ \eqref{eq:splitpsi} to decompose each field into the potential from a \emph{finite} number of charges/images and a harmonic correction. In order to solve the correction problems analytically, we require that the correction $\cor{\psi}_b$ be harmonic on $(-\infty,0]$, and that the correction $\cor{\psi}_t$ be harmonic on $[H,\infty)$. This means that the intermediate potential $\psi_b^*$ must include any charges whose support extends below $z=0$ (those centered $H_E$ or closer to the bottom wall, like the bottom-most blue charge in the right column of Fig.\ \ref{fig:ffsolve}), and the intermediate potential $\psi_t^*$ must include any charges whose support extends above $z=H$ (those centered $H_E$ or closer to the top wall). To minimize the number of doubly-periodic Poisson solves, we combine these charges into an intermediate potential $\psi_o^*$ which includes any charges $\olap{C}$ centered within $2H_E$ of the top and bottom boundaries (blue charges in the right column of Fig.\ \ref{fig:ffsolve}), 
\begin{equation}
\label{eq:phiost}
    \epsilon \Delta \psi_o^* = \olap{\farf}. 
\end{equation}
This Poisson equation can again be solved using the method of Section\ \ref{sec:p1smooth}. Since all of the charges included are centered within the slab and have support $H_E$ in all directions, we could use the domain $z \in [-H_E,H+H_E]$. The algorithm is simpler and more efficient, however, if we use the same domain and grid as the Poisson solve\ \eqref{eq:phiistar}, $z \in [-3H_E,H+3H_E]$. If we assign the charges $\olap{C}$ the proper image strength (depending on whether we seek the solution above or below the slab), we obtain the representation of the intermediate potentials
\begin{gather}
\label{eq:phibst}
\psi^*_b =  \frac{2\epsilon}{\epsilon_b+\epsilon}\psi_o^*, \\
\label{eq:phitst}
\psi^*_t =  \frac{2\epsilon}{\epsilon_t+\epsilon}\psi_o^*. 
\end{gather}

\subsubsection{Correction potentials \label{sec:farharm}}
In our initial solve\ \eqref{eq:phiistar} for $\psi_i^*$, we did not include the images $\nolap{C}^\text{(img)}$, which are shown as dashed red circles in Fig.\ \ref{fig:ffsolve}. Subtracting\ \eqref{eq:phiistar} from the original equation\ \eqref{eq:phii} for $\psi_i$, we get a Poisson equation for the correction $\cor{\psi}_i=\psi_i-\psi_i^*$ which involves these images,  
\begin{equation}
\epsilon \Delta \cor{\psi}_i = -\nolap{\farf}^\text{(img)}.
\end{equation}
Since $\nolap{\farf}^\text{(img)}$ is compactly supported outside of $z \in [-H_E,H+H_E]$ (no overlap of the dashed empty circles with the gray region in Fig.\ \ref{fig:ffsolve}), we have the harmonic problem for the interior correction 
\begin{equation}
\label{eq:harmi}
\epsilon \Delta \cor{\psi}_i = 0  \text{ for } z \in [-H_E,H+H_E].
\end{equation}

For the exterior corrections, the splitting of the r.h.s.\ $\farf$ in\ \eqref{eq:fsp} and definition of $\psi_b$ in\ \eqref{eq:psib} imply that the correction $\cor{\psi}_b$ solves the Poisson equation
\begin{gather}
\label{eq:RHSb}
\epsilon \Delta \cor{\psi}_b = -\frac{2\epsilon}{\epsilon_b+\epsilon} \left(\nolap{\farf}+\img{\farf}_{z>H}\right).
\end{gather}
Since the r.h.s.\ of\ \eqref{eq:RHSb}, which comprises the dashed open circles positioned above $z=0$ in the right column of Fig.\ \ref{fig:ffsolve}, is compactly supported far away from the bottom wall (all charges and images are centered at least $2H_E$ from the bottom wall), we obtain the Laplace equation
\begin{equation}
\label{eq:harmb}
\epsilon \Delta \cor{\psi}_b= 0  \text{ on }(-\infty,0].
\end{equation}
The procedure for $\cor{\psi}_t$ is identical. We observe that the empty charges and images positioned below $z=H$ in the right column of Fig.\ \ref{fig:ffsolve} are far from the top wall of the slab. The correction field that they generate is therefore harmonic above $z=H$, 
\begin{equation}
\label{eq:harmt}
\epsilon \Delta \cor{\psi}_t= 0  \text{ on }[H,\infty).
\end{equation}
We now have the three Laplace equations\ \eqref{eq:harmi},\ \eqref{eq:harmb}, and\ \eqref{eq:harmt} (on different domains) for $\cor{\psi}_i$, $\cor{\psi}_b$, and $\cor{\psi}_t$, just as in Section\ \ref{sec:harmsolve}. 

The final step is to impose boundary conditions for the correction solve such that the total fields $\psi_i$, $\psi_o$, and $\psi_t$ satisfy the boundary conditions\ \eqref{eq:psiBC1}$-$\eqref{eq:psiBC4}. If we use the mismatches
\begin{gather}
\label{eq:mm1}
m^{(b)}_\phi(x,y)=\psi_i^*(x,y,z=0)-\psi_b^*(x,y,z=0),\\
\label{eq:mm2}
m^{(b)}_E(x,y) = \epsilon \frac{\partial \psi_i^*}{\partial z}(x,y,z=0)-\epsilon_b \frac{\partial \psi_b^*}{\partial z}(x,y,z=0)+\sigma_b(x,y),\\
m^{(t)}_\phi(x,y)=\psi_i^*(x,y,z=H)-\psi_t^*(x,y,z=H),\\
\label{eq:mm4}
m^{(t)}_E(x,y)=\epsilon \frac{\partial \psi_i^*}{\partial z}(x,y,z=H)-\epsilon_t \frac{\partial \psi_t^*}{\partial z}(x,y,z=H)-\sigma_t(x,y),
\end{gather}
as boundary conditions for the harmonic solve described in Section\ \ref{sec:harmsolve}, we obtain $\cor{\psi}_i$ on $[-H_E,H+H_E]$ and, if desired, $\cor{\psi}_t$ on $z > H$ and $\cor{\psi}_b$ on $z<0$, for all modes except $\kpar=0$. We could combine the corrections $\cor{\psi}_{i/t/b}$ with the intermediate potentials $\psi^*_{i/t/b}$ to obtain $\psi$ everywhere, although we only need $\psi_i$ on $[-H_E,H+H_E]$ to evaluate the energy and forces on the charges. In any case, the $\kpar=0$ mode must be handled separately as discussed in the next section.

\subsubsection{Corrections for $\kpar=0$ \label{sec:farcor}}
In a natural extension of the method presented in Section\ \ref{sec:combsmooth}, we obtain the solution for $\kpar=0$ by considering the total solutions $\psi_{i/t/b}=\psi_{i/t/b}^*+\cor{\psi}_{i/t/b}$, rather than considering the intermediate and correction potentials separately. In $xy$ Fourier space, we have the solutions for the $\kparv=\V{0}$ mode
\begin{gather}
\label{eq:psiicor}
\four{\psi}_i (\kparv=\V{0},z)= \four{\psi}_i^*(\V{0},z)+A_i(\V{0}) z + B_i(\V{0}),\\
\label{eq:psibcor}
     \four{\psi}_b (\kparv=\V{0},z)=  \four{\psi}_b^*(\V{0},z)+A_b(\V{0})z + B_b(\V{0}),\\
     \label{eq:psitcor}
    \four{\psi}_t (\kparv=\V{0},z)=  \four{\psi}_t^*(\V{0},z)+A_t(\V{0})z+B_t(\V{0})
\end{gather}
By $\four{\psi}_{i/t/b}^*(\V{0},z)$, we mean the $\kpar=0$ components of the solutions of the Poisson equations\ \eqref{eq:phiistar} and\ \eqref{eq:phiost}, 
\begin{gather}
\label{eq:phii0bvp}
    \four{\psi}_i^*(\V{0},z)=-\frac{1}{\epsilon  }\int_{-3H_E}^z \, dz' \int_{-3H_E}^{z'} \, dz'' \left( \olap{\four{\farf}}(\V{0},z'')+\nolap{\four{\farf}}(\V{0},z'') +\olap{\four{\farf}}^\text{(img)}(\V{0},z'')\right),\\
    \label{eq:phio0bvp}
     \four{\psi}_o^*(\V{0},z)=-\frac{1}{\epsilon}\int_{-3H_E}^z \, dz' \int_{-3H_E}^{z'} \, dz'' \;  \olap{\four{\farf}}(\V{0},z'').
\end{gather}
Equation\ \eqref{eq:phio0bvp} defines $\four{\psi}_o^*(\V{0},z)$, from which  $\four{\psi}_b^*(\V{0},z)$ and $\four{\psi}_t^*(\V{0},z)$ can be obtained by\ \eqref{eq:phibst} and\ \eqref{eq:phitst}. Note that the formulas\ \eqref{eq:phii0bvp} and\ \eqref{eq:phio0bvp} imply a set of boundary conditions for the corresponding BVPs\ \eqref{eq:phiistar} and\ \eqref{eq:phiost} for $\kpar=0$; these BCs are arbitrary since $\four{\psi}_{i/o}^*(\V{0},z)$ is only defined up to a linear function. 

Once the solutions $\four{\psi}_{i/b/t}^*(\V{0},z)$ are obtained, we can solve for the coefficients $A_b(\V{0})$ and $A_t(\V{0})$ by imposing the decay boundary condition for $\four{\psi}_b$ in\ \eqref{eq:psibcor} and $\four{\psi}_t$ in\ \eqref{eq:psitcor} to obtain
\begin{equation}
\label{eq:CE}
A_b(\V{0}) = -\frac{\partial \four{\psi}_b^*}{\partial z}\left(\V{0},z=-3H_E\right), \qquad A_t(\V{0}) = -\frac{\partial \four{\psi}_t^*}{\partial z}\left(\V{0},z=H+3H_E\right). 
\end{equation}
The electric displacement boundary conditions\ \eqref{eq:psiBC2} and \eqref{eq:psiBC4} now give $A_i(\V{0})$, 
\begin{align}
\label{eq:ai1}
    \epsilon A_i(\V{0}) & =  \epsilon_b \left( \frac{\partial \four{\psi}_b^*}{\partial z}\left(\V{0},0\right)+ A_b(\V{0})\right) -\epsilon \frac{\partial \four{\psi}_i^*}{\partial z}\left(\V{0},0\right)-\four{\sigma}_b(\V{0}) \\ 
    \label{eq:ai2}
     & = \epsilon_t \left( \frac{\partial \four{\psi}_t^*}{\partial z}\left(\V{0},H\right)+ A_t(\V{0})\right) -\epsilon \frac{\partial \four{\psi}_i^*}{\partial z}\left(\V{0},H\right)+\four{\sigma}_t(\V{0}).
\end{align}
For electroneutral slabs, the right hand sides of\ \eqref{eq:ai1} and\ \eqref{eq:ai2} must be the same in continuum.\footnote{This can be verified by writing out both sides in terms of the charge densities, substituting the solutions\ \eqref{eq:phii0bvp}$-$\eqref{eq:CE}, then using the identity $\int_0^H \olap{\farf}^\text{(img)} \, dz = (\epsilon-\epsilon_b)/(\epsilon_b+\epsilon)\int_{-3H_E}^0 \olap{\farf} \, dz$ for far-field images that overlap the bottom wall. Using a similar identity for images overlapping the top wall and the electroneutrality condition\ \eqref{eq:elecneu} then gives the equality\ \eqref{eq:ai1}$=$\eqref{eq:ai2}. }
Finally, as in Section\ \ref{sec:combsmooth}, we have the continuity conditions\ \eqref{eq:psiBC1} and\ \eqref{eq:psibc3} which give
\begin{gather}
\label{eq:Bb}
    B_b(\V{0})=\four{\psi}_i^*(\V{0},0)+B_i(\V{0})-\four{\psi}_b^*(\V{0},0),\\
\label{eq:Bt}
B_t(\V{0})=\four{\psi}_i^*(\V{0},H)-\four{\psi}_t^*(\V{0},H)+\left( A_i(\V{0})-A_t(\V{0}) \right)H+B_i(\V{0}).
\end{gather}
In Section\ \ref{sec:ewslabalg}, we describe how to obtain $B_i(\V{0})$ in real space according to\ \eqref{eq:phi00}. The values of $B_b(\V{0})$ and $B_t(\V{0})$ then follow from\ \eqref{eq:Bb} and\ \eqref{eq:Bt}. 

\section{Numerical method and algorithm \label{sec:numerics}}
In this section, we discuss the implementation of the algorithm we described in Section\ \ref{sec:ourstuff}. Here we focus on getting $3-4$ digits of accuracy, as this is more than sufficient for Brownian dynamics simulations, but it is straightforward to achieve higher accuracies if desired. We begin in Section\ \ref{sec:cutgrid} with precomputations, in particular our initialization of the grid and cutoff distances for near and far field truncation. We move on to main computations in Section\ \ref{sec:alg}, where we first discuss our solver for general doubly-periodic problems and then specialize to the Ewald splitting of Section\ \ref{sec:ourstuff}. The number of images we include in the near field and far field calculations have important implications for the allowed values of the Ewald parameter $\xi$, as we explain in Section\ \ref{sec:xi}. 

\subsection{Cutoffs and grid spacing \label{sec:cutgrid}}
The splitting parameter $\xi$ determines the truncation distances for the near field kernel and far-field charge width. Beginning with the near field, we will take the near field cutoff distance $\rnf$ to be the minimum value such that (see\ \eqref{eq:Gnear} and\ \eqref{eq:Enear})
\begin{equation}
\label{eq:rnf}
\abs{\frac{\partial \near{G}(\rnf;g_w,\xi)/\partial r }{\partial \near{G}(\rnf;g_w,0)/\partial r }} < \delta,
\end{equation}
where $\delta =10^{-4}$ for $\sim 4$ digits of accuracy and $\delta=5 \times 10^{-4}$ for $\sim 3$ digits. Recall that the average near field kernels\ \eqref{eq:GnearK} and\ \eqref{eq:EnearK} are averages of the pointwise near field kernels over a Gaussian with support $n_\sigma g_w$. As such, we truncate the kernels\ \eqref{eq:GnearK} and\ \eqref{eq:EnearK} for $r> \rcut = \rnf+n_\sigma g_w$. 

We use a Fourier-Chebyshev grid to solve smooth doubly periodic problems following the approach of Section\ \ref{sec:smoothdp}. For the far field discretization, we use FFTs in the $x$ and $y$ directions, so we take the (FFT-friendly) number of grid points $N_{x/y}$ as an input, which we combine with $g_w$ to calculate the total far-field charge width $g_t$ using\ \eqref{eq:He} and splitting parameter $\xi$ using\ \eqref{eq:gt}.  We consider two different uniform grid spacings in $x$ and $y$, with $h_{xy} =g_t/1.2$ giving $\sim 3$ digits of accuracy and $h_{xy}=g_t/1.4$ giving $\sim 4$ digits.  In the examples studied here, $h_x=L_x/N_x=h_y=L_y/N_y=h_{xy}$. 

Although the Gaussian\ \eqref{eq:spreadkernelEw} is technically nonzero everywhere on the grid, we truncate it at a finite number of grid cells $n_g$ for efficiency in spreading the charge density to the far field grid (because the $x$ and $y$ grid is uniformly spaced, this assumption allows for the possible optimization of fast Gaussian gridding \cite{NUFFT} in the $x$ and $y$ directions). This gives the total radius of the Gaussian spreading kernel as
\begin{equation}
\label{eq:He}
H_E=\frac{n_g}{2}h_{xy} = n_\sigma g_t,
\end{equation}
where $n_\sigma$ is the (non-integer) number of standard deviations over which the Gaussian is supported, and the radius of support of the Gaussian charge is $H_E$. For $\sim 3$ digits of accuracy, we will use a Gaussian support $n_g=10$, which gives a Gaussian truncated at $n_\sigma =5/1.2 \approx 4.2$ standard deviations. For 4 digits, we take $n_g=12$ for $n_\sigma=6/1.4 \approx 4.3$ standard deviations. \rev{Since our method is formally spectrally accurate prior to Gaussian truncation, in principle, arbitrarily high accuracy (up to roundoff) can efficiently be achieved by using larger Gaussian truncation widths and more grid cells in the support of the Gaussian kernels.} 

Once the $xy$ Fourier grid is chosen, the number of Chebyshev grid points for the $z$ grid can be chosen based on accuracy considerations. Specifically, we require the spacing at the middle of the Chebyshev $z$ grid (where it is the coarsest) to be comparable to that of the Fourier $xy$ grid. In a Chebyshev grid, the maximum spacing in the middle of the grid is at most $\pi/2$ times the average spacing, \rev{so our Chebyshev grid requires about $\pi/2$ as many points as a uniform grid with the same maximum spacing. Specifically,} if $\widebar{H}$ is the total height of the domain in a doubly periodic solve as described in Section\ \ref{sec:smoothdp}, we set the number of Chebyshev points to
\begin{equation}
N_z = \left\lceil \frac{\pi \widebar{H}}{2h_{xy}}\right\rceil,
\end{equation} 
where the integer rounding should be chosen for FFT optimality.\footnote{The complex FFTs which we employ here to transform from values on the grid to Chebyshev coefficients have $2N_z-2$ points in the $z$ direction, so rounding should be such that $2N_z-2$ is favorable for the FFT.} Since the $z$ (Chebyshev) grid is irregularly spaced, fast Gaussian gridding is not an option, and we must simply assign the kernel\ \eqref{eq:spreadkernelEw} a value of zero at any point farther than $ H_E$ in the $z$ direction from the charge center. \rev{The fixed truncation distance implies that we must evaluate the kernel\ \eqref{eq:spreadkernelEw} at more grid points near the domain boundaries, since the Chebyshev grid is finer there. }

\subsection{Algorithm \label{sec:alg}}
The far field solver described in Section\ \ref{sec:farfield} requires solving two doubly periodic Poisson equations\ \eqref{eq:phiistar} and\ \eqref{eq:phiost}, and so we begin this section by describing an algorithm to solve them using the method developed in Section\ \ref{sec:p1smooth}. We then summarize the complete algorithm for slabs. A GPU implementation of our algorithm is available freely at github, see \url{https://github.com/stochasticHydroTools/DPPoissonTests/} for instructions and examples.

\subsubsection{Boundary value solver \label{sec:dpalg}}
Our doubly periodic Poisson solver is based on transforming the charge density $f$ into its Fourier/Chebyshev representation on the grid, and consists of the following steps:
\begin{enumerate}
\item Compute the charge density $f(x,y,z)$ on the Fourier-Chebyshev grid. 

\item Take the 3D FFT of the charge density to obtain the Fourier and Cheyshev coefficients on the grid. We refer to this as a fast Fourier-Chebyshev transform (FFCT). See \cite[c.~8]{trefethen2000spectral} for a description of how to obtain Chebyshev coefficients using the FFT. 

\item Use the Chebyshev boundary value solver \cite{greengard1991spectral} described in Appendix\ \ref{sec:bvps} to solve the BVPs\ \eqref{eq:BVP} for each wave number $\kparv$ using a well-conditioned integral formulation. This gives the Fourier-Chebyshev coefficients of the potential $\phi$. Note that this step is trivially parallelizable since each mode $\kparv$ is handled independently of others.

\item For the electric field, compute derivatives of $\phi$ on the grid by Fourier differentiation in $x$ and $y$ (i.e., by multiplying $\four{\phi}$ by $ik_x$ or $ik_y$), or differentiating the Chebyshev series in $z$. This gives the Fourier-Chebyshev coefficients of $\nabla \phi$ on the grid. 
\item Perform a 3D inverse fast Fourier-Chebyshev transform (IFFCT) to obtain $\phi$ and $\V{E}=-\nabla \phi$ on the $(x,y,z)$ grid. 
\end{enumerate}
Note that if only energy or only forces are required then some of these steps can be skipped. 

\subsubsection{Ewald splitting for slabs \label{sec:ewslabalg}}
We now detail our Ewald splitting algorithm for computing $\widebar{\phi}$ and $\widebar{\V{E}}$ at the centers $\V{z}_k$ of all charges in the slab geometry. The near field algorithm is straightforward. From $\xi$, we determine $\rnf$ by solving\ \eqref{eq:rnf} and then, for each charge and dielectric boundary, we construct one image and compute the sum\ \eqref{eq:avgGn} and/or\ \eqref{eq:EnearK} for each charge. This sum includes only pairs of charges closer than $\rcut$ apart, and therefore the cost of the near field algorithm is linear in the number of particles.  

The far field algorithm is more complex and worth listing in steps, with a graphic representation given in Fig.\ \ref{fig:ffsolve}. Because image charges centered up to $2H_E$ away from the slab can be included in some of the solves and each smeared (far-field) charge has a support of $H_E$ in all directions, we initialize a single second-kind Chebyshev grid on $z \in [-3H_E,H+3H_E]$ and a Fourier grid on $[0,L_x] \times [0,L_y]$. We then perform the steps:

\begin{enumerate}
\item Compute $H_E$ using\ \eqref{eq:He}. Then separate the charges into two groups: charges $\olap{C}$ with $z < 2H_E$ or $z>H-2H_E$ (solid blue circles in Fig.\ \ref{fig:ffsolve}), and the rest of the charges $\nolap{C}$ (solid green circles in Fig.\ \ref{fig:ffsolve}). Then compute the positions and strengths of the necessary images $\olap{C}^\text{(img)}$ (solid red circles in Fig.\ \ref{fig:ffsolve}) using\ \eqref{eq:ims}. 

\item Construct the intermediate potential $\psi_o^* $ outside the slab by solving the doubly periodic problem\ \eqref{eq:phiost} as described in Section\ \ref{sec:smoothdp} for $z \in [-3H_E,H+3H_E]$. Specifically, first spread the charge $q_i$ for $i \in \olap{C}$ onto the grid using the kernel $\widebar{S}$ defined in\ \eqref{eq:spreadkernelEw}. Then apply steps $2-4$ in Section\ \ref{sec:dpalg} to obtain $\four{\psi}_o^* $ on the Fourier-Chebyshev grid. For simplicity, use the BVP solver to also obtain the part of the $\kpar=0$ mode\ \eqref{eq:phio0bvp} by solving the corresponding BVP\ \eqref{eq:phiost} for $\kpar=0$ with homogeneous Dirichlet BCs at $z=-3H_E$ and $z=H+3H_E$. The intermediate potentials $\four{\psi}_b^* $ and $\four{\psi}_t^* $ in Fourier/Chebyshev space can trivially be obtained by multiplying $\four{\psi}_o^* $ by the coefficients in\ \eqref{eq:phibst} and\ \eqref{eq:phitst}. 

\item Construct the intermediate potential $\psi_i^*$ by solving the doubly periodic problem\ \eqref{eq:phiistar} for $z \in [-3H_E,H+3H_E]$. To do this, spread the charges $q_i$ for $i \in \nolap{C}$ and images $q_j$ for $j \in \olap{C}^\text{(img)}$ to the grid and add the result to the spreading for $\olap{C}$ already computed in step 2. Then apply steps $2-4$ in Section\ \ref{sec:dpalg} to obtain $\four{\psi}_i^* $ on the Fourier-Chebyshev grid. For simplicity, use the BVP solver to obtain the part of the $\kpar=0$ mode\ \eqref{eq:phii0bvp} by solving the corresponding $\kpar=0$ BVP\ \eqref{eq:phiistar} with homogeneous Dirichlet BCs at $z=-3H_E$ and $z=H+3H_E$.

\item Calculate the mismatches $m_\phi(\kparv \neq \V{0})$ and $m_E(\kparv \neq \V{0})$ given in\ \eqref{eq:mm1}$-$\eqref{eq:mm4} using Chebyshev differentiation in $z$. Use these mismatches as boundary conditions to obtain the harmonic corrections $\cor{\psi}_i$, $\cor{\psi}_b$, and $\cor{\psi}_t$ analytically as outlined in Section\ \ref{sec:harmsolve}. Evaluate the solution for each $\kparv$ and for each Chebyshev grid point $z_c$ with $-H_E \leq z_c \leq H+H_E$, and set $\cor{\psi}_i(\kparv,z)=0$ outside of this $z$ range to avoid over/underflow errors. Because of the ill-conditioning of the correction solve for large $\kpar$, set $\cor{\psi}_i(\kparv,z)$ to zero for $\kpar > k_\text{max} = \pi/h_\text{xy}$. Finally, perform $N_x N_y$ independent 1D Chebyshev transforms (FFTs) in $z$ to obtain the Fourier-Chebyshev representation of $\cor{\four{\psi}}_i$.\footnote{The transforms from $z$ space to Chebyshev space are done so that $\cor{\four{\psi}}_i$ can be combined with $\four{\psi}_i^*$ directly in Fourier-Chebyshev space, and then a single 3D IFFCT performed to obtain $\psi_i$. An alternative but less efficient sequence is to transform $\cor{\four{\psi}}_i$ directly from Fourier-$z$ space to real space via a (parallel) series of 2D IFFTs in the $xy$ plane, then add to the result from a 3D IFFCT on $\four{\psi}_i^*$.} Note that this step can be trivially parallelized since each mode $\kparv$ is handled independently. 

\item  To correct the $\kpar=0$ mode, modify the solution for $\psi_i^*(\kparv=\V{0},z)$ already obtained in step 4 by adding the linear mode $A_i(\V{0})z$ as discussed in Section\ \ref{sec:farcor}. In the discrete setting, the equations for $A_i(\V{0})$,\ \eqref{eq:ai1} and\ \eqref{eq:ai2}, give the same value for $A_i(\V{0})$ to only about 3 relative digits, so set $A_i(\V{0})$ to be the mean result of the two.

\item Set the far field values in Fourier-Chebyshev space to $\four{\psi}_i = \four{\psi}_i^*+\cor{\four{\psi}}_i$. Then perform a 3D IFFCT to obtain $\four{\psi}_i$ on the grid. 

\item Interpolate (average) $\psi_i$ at the charge centers using the kernel $\widebar{S}$ given in\ \eqref{eq:spreadkernelEw} to obtain the average far field potential $\far{\widebar{\phi}}$ at the charge centers.
\end{enumerate}

The calculation of the electric field follows a similar procedure. Since the potential on the grid $\psi_i =\psi_i^*+\cor{\psi}_i$ has two parts, we obtain the electric field inside the slab by differentiating each. We first differentiate $\psi_i^*$ on the grid using Fourier multiplication and Chebyshev differentiation. We then differentiate $\cor{\psi}_i$ analytically in $z$ and with Fourier multiplication for $x$ and $y$. Adding these two pieces together gives $-\nabla \psi_i$, which we then interpolate back onto the charges to obtain the average electric field $\far{\widebar{\V{E}}}$ on each charge. Note that another way to obtain the electric field is to interpolate the potential $\psi_i$ with the derivative of the interpolation kernel $\widebar{S}$, but this enlarges the width of the interpolation kernel, thereby requiring interpolation from more grid cells to obtain the same accuracy.

When the energy is needed in addition to the forces, we add a constant to $\phi$ (represented by $B_i(\V{0})$) to set the pointwise potential $\phi(\V{0})=0$ and remove an arbitrary shift in the energy. To do this, we sum the pointwise near field kernel at $\V{x}=\V{0}$ using the kernel\ \eqref{eq:Gnearsum} (with truncation at $\rnf$ and one set of images above and below the slab) to obtain $\near{\phi}(\V{0})$. For the far field, we interpolate the potential $\psi_i$ with the kernel $\gamma^{1/2}$ at $\V{x}=\V{0}$ to obtain the pointwise far field potential $\far{\phi}(\V{0})$. We then set $B_i(\V{0})$ in\ \eqref{eq:psiicor} so that the total potential $\near{\phi}(\V{0})+\far{\phi}(\V{0})=0$. 

\subsection{Constraints on the splitting parameter \label{sec:xi}}
Our decision to include only the first set of images in the near field and far field problems leads to restrictions on the Ewald parameter $\xi$. In both of these problems, we make the assumption that images above the top wall cannot interact with points at or below the bottom wall, and vice versa. This constrains the truncation distances $\rnf$ and $H_E$, and therefore $\xi$ itself.  

\subsubsection{Relationship between $\xi$ and $g_w$}
In order for there to be a reduction in the far field problem grid size significant enough to justify Ewald splitting in the first place, the smearing\ \eqref{eq:fmod} must increase the total width of the Gaussian charge cloud by a substantial amount. Here we take substantial to mean that the smeared Gaussian width $g_t$ is at least twice as large as the original width $g_w$, so that $g_t \geq 2g_w$. Using\ \eqref{eq:gt} to solve for $\xi$, we get the bound
\begin{equation}
\label{eq:xieff}
g_w \sqrt{12} \leq \frac{1}{\xi} \Leftrightarrow \xi \leq \frac{1}{g_w\sqrt{12}}.
\end{equation}

\subsubsection{Near field constraints \label{sec:nearconstraint}}
We recall the near field constraint that images above the slab cannot interact with the bottom wall, $\rnf < n_\text{img}H+h$, where $h \geq 4g_w$ is the minimum distance between a charge and one of the dielectric boundaries. This assumption gives a constraint on $\xi$ as well. To estimate the bound, we set $\rnf = \alpha/\xi$ and use the upper bound on $\xi$ from\ \eqref{eq:xieff} to obtain the equation
\begin{equation}
\label{eq:alphest}
\delta = \abs{\frac{\partial \near{G}(\alpha/\xi;1/(\sqrt{12}\xi),\xi)/\partial r }{\partial \near{G}(\alpha/\xi;1/(\sqrt{12}\xi),0)/\partial r }}, 
\end{equation}
which is independent of $\xi$ and can be solved numerically for $\alpha$ to obtain a bound
\begin{equation}
\label{eq:nearapprox}
n_\text{img}H+h > \frac{\alpha}{\xi} > \rnf . 
\end{equation}
The actual value of $\alpha$ depends on the relationship between the smeared Gaussian width $g_t$ and the original Gaussian width $g_w$. For $\delta=10^{-4}$, $\alpha \approx 3.5$ when $g_t/g_w \geq 2$, while $\alpha=3.25$ for $g_t \gg g_w$, which gives the range $\alpha \in [3.25,3.50]$. For $\delta = 5 \times 10^{-4}$, we get $\alpha \in [2.98,3.22]$.

\subsubsection{Far field constraints}
We next consider the far field constraints on $\xi$. In the far field solves, we never include images of images explicitly, always assuming they can be included in the correction solve because they are 
a distance $H+h$ from the slab boundaries. For this to be the case, the images of images cannot overlap the region $z \in [-H_E, H+H_E]$, meaning they must be at least $2H_E \geq 8g_t$ away from the slab, so that in total we have the constraint $H+h > 2H_E$. 

To estimate $H_E$ in terms of $\xi$, we first estimate $g_t$. For Ewald splitting to actually reduce the computational complexity, we recall from\ \eqref{eq:xieff} that we require $g_w \leq 1/(\xi \sqrt{12})$. This gives an upper bound on the total Gaussian width $g_t$, 
\begin{equation}
g_t = \sqrt{\frac{1}{4\xi^2}+g_w^2} \leq \frac{1}{\xi \sqrt{3}}. 
\end{equation}
In the case $g_w \rightarrow 0$, we get $g_t = 1/(2\xi)$. Since images of images are at least $H+h$ away from either wall, we need 
\begin{equation}
\label{eq:farconstraint}
2H_E =2n_\sigma g_t \leq \frac{2n_\sigma}{\xi \sqrt{3}} < H+h \Rightarrow \frac{1}{\xi} < \frac{\sqrt{3}(H+h)}{2n_\sigma}, 
\end{equation} 
or $1/\xi < (H+h)/n_\sigma$ in the case $g_w \rightarrow 0$. 

The constraints for Ewald splitting are summarized in Table\ \ref{tab:ewconstraints}. We see that, when we do not include images of images (i.e., $n_\text{img}=1$), the far field and near field constraint on $\xi$ are almost identical, with the far field constraint being more restrictive. Thus the far field calculation restricts the Ewald parameter so much that including more images in the near field is not necessary.

\begin{table}
\centering
\begin{tabular}{c|c|c}
Tolerance $\delta$ &  $10^{-4}$ & $5 \times 10^{-4}$\\
\hline
Grid size $h_{xy}$ & $g_t/1.4$ & $g_t/1.2$\\
Gaussian support $n_g$ & 12 & 10 \\
Gaussian truncation $n_\sigma$ & 4.29 & 4.17\\
Near field constraint $ 1/\xi < $ & $(n_\text{img} H+h)/[3.25,3.50]$ &$(n_\text{img} H+h)/[2.98,3.22]$\\
Far field constraint $1/\xi < $ & $(H+h)/[4.29,4.95]$ & $(H+h)/[4.17,4.81]$\\
\end{tabular}
\caption{Summary of constraints for Ewald splitting for $\sim 4$ (left column) and $\sim 3$ (right column) digits of accuracy. For each tolerance, we first give the grid spacing $h_{xy}$, number of grid cells $n_g$ in the support of the far-field Gaussian, and the corresponding number of Gaussian standard deviations $n_\sigma$ in the truncated support. We then give bounds on $1/\xi$ from the near-field constraint\ \eqref{eq:nearapprox} and far-field constraint\ \eqref{eq:farconstraint}. Lower ends of the intervals (weaker constraints) are for the case when $g_w \rightarrow 0$ (point charges), and upper ends are for the case when $g_t=2g_w$. }
\label{tab:ewconstraints}
\end{table}


\section{Numerical tests \label{sec:tests}}
In this section, we validate our algorithm for dielectric slabs. We begin in Section\ \ref{sec:fscomp} with a simple verification: four charges inside of a slab unbounded in the lateral directions. In Section\ \ref{sec:ewind}, we verify that our answer is independent of the Ewald splitting parameter $\xi$ by generating a reference result \emph{without} Ewald splitting ($\xi \rightarrow \infty$). We conclude in Section\ \ref{sec:fge} by verifying that the force\ \eqref{eq:fqE} is indeed the gradient of the energy\ \eqref{eq:Ufinal}. Throughout this section, we use a tolerance $\delta=10^{-4}$ (which gives the finer grid setting $n_g=12$, $h_{xy} = g_t/1.4$) for tests with 10 or fewer charges (Sections\ \ref{sec:fscomp} and\ \ref{sec:fge}), and a tolerance $\delta=5 \times 10^{-4}$ ($n_g=10$, $h_{xy}=g_t/1.2$) for tests with more than 10 charges (Section\ \ref{sec:ewind}).

\subsection{Comparison with the image construction in free space \label{sec:fscomp}}
We first verify that our spectral Ewald method agrees with the free space solution as the periodic $xy$ length $L_x=L_y=L \rightarrow \infty$. The solution in free space can be derived via the image construction with an infinite series of images (see e.g. \cite{tinkle2001image}). 

We set up our test with four charges positioned randomly in the $xy$ plane and with $z$ locations 0.01, 0.30, 1.00, 1.95 (see Table \ref{tab:ewfreeslab}). We set $\epsilon=1$, $\epsilon_t=1/5$, $\epsilon_b=1/2$ so that the free space image sums converge rapidly (we use a total of 400 image charges so that the partial sums are accurate to machine precision), $H=2$, and $g_w=10^{-2}$. We let $\xi \approx 3$ ($g_t=0.166$, $H_E=0.71$, $\rnf=1.08$), solve the problem for $L=28$ ($236 \times 236 \times 84$ grid) and $L=32$ ($270 \times 270 \times 84$ grid), extrapolate the result to $L=\infty$ (based on the fact that the finite-size correction decays like $1/L$), and denote the resulting electric field by $\widebar{\bm{E}}_E$. We show the normalized difference $\widebar{\bm{E}}_E-\widebar{\bm{E}}_F$, where $\widebar{\bm{E}}_F$ is the free space result, in Table \ref{tab:ewfreeslab}, which shows that we obtain about 5 digits of accuracy for all four charges.

\begin{table}
\centering
\begin{tabular}{c|c|c}
$\V{z}_i$ & $q_i$ &$\left(\widebar{\bm{E}}_E(\V{z}_i)-\widebar{\bm{E}}_{F}(\V{z}_i)\right)/\text{mean} \left(\norm{\widebar{\V{E}}_F(\V{z}_k)}\right)$ \\
\hline
(-0.17,-0.71,0.01)& 1 & (-9.2e-6, -1.9e-6, 2.4e-6)\\
(0.44,-0.82,0.30)& -1 & (-9.4e-6, -1.9e-6, 3.0e-6)\\
(-1.00,-0.63,1.00)& 1 & (-9.8e-6, -1.5e-6, 3.1e-6)\\
(-0.40,-0.31,1.95)& -1 & (-9.5e-6, -1.6e-6, 2.5e-6)
\end{tabular}
\caption{Relative errors in the extrapolated values of the electric field from the Ewald method as $L_x=L_y \rightarrow \infty$ vs.\ the image method for free space, $\widebar{\bm{E}}_E-\widebar{\bm{E}}_{F}$, for four charges positioned at $\V{z}_i$ with strength $q_i$.  }
\label{tab:ewfreeslab}
\end{table}

\subsection{Independence from splitting parameter \label{sec:ewind}}
We next verify that our results are independent of the splitting parameter $\xi$. To do this, we compare the results for the spectral Ewald method with various $\xi$ to reference results obtained without Ewald splitting ($\xi \rightarrow \infty$). We generate the reference results to 7 digits of accuracy using the method of Section\ \ref{sec:smoothdp} with a fine grid. 

We consider 100 random charges of unit strength and alternating sign positioned on $[-1,1] \times [-1,1] \times [4.5g_w, H-4.5g_w]$ and set $g_w=0.025$ so that it is not too expensive to generate the reference solution to high accuracy. For the jump conditions, we set $\epsilon=1$, $\epsilon_t=1/50$, and $\epsilon_b=1/20$. We take $H=0.75$, so that, for some choices of the Ewald parameter, charges in the center of the slab have (smeared far-field) support that intersects both $z=0$ and $z=H$. 

We choose a set of Ewald parameters in the range $\xi=4$ to $\xi=26$. The lower bound on $\xi$ is chosen so that each point interacts only with its nearest image to 4 digits in the near field calculation (this includes periodic images as well as images in the $z$ direction). The upper bound on $\xi$ is the point where the grid is almost as large as that used for the calculation without Ewald splitting (i.e., $\xi$ is so large that the Ewald splitting is useless). Table \ref{tab:ewnoparam} summarizes the values of the Ewald parameter we use and the resulting values of $H_E$ and $\rnf$.

\begin{figure}
\centering     
\includegraphics[width=0.45\textwidth]{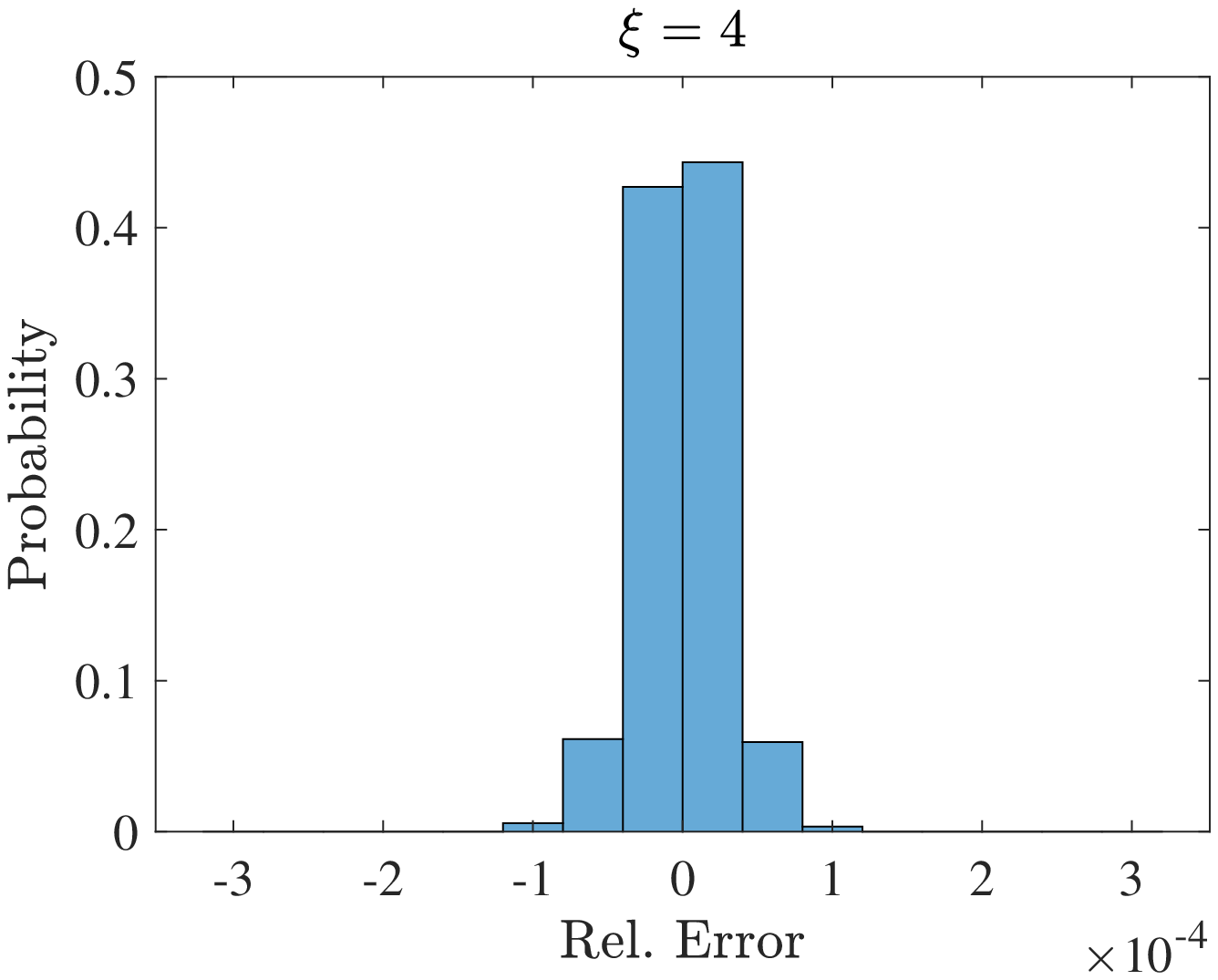}
\includegraphics[width=0.45\textwidth]{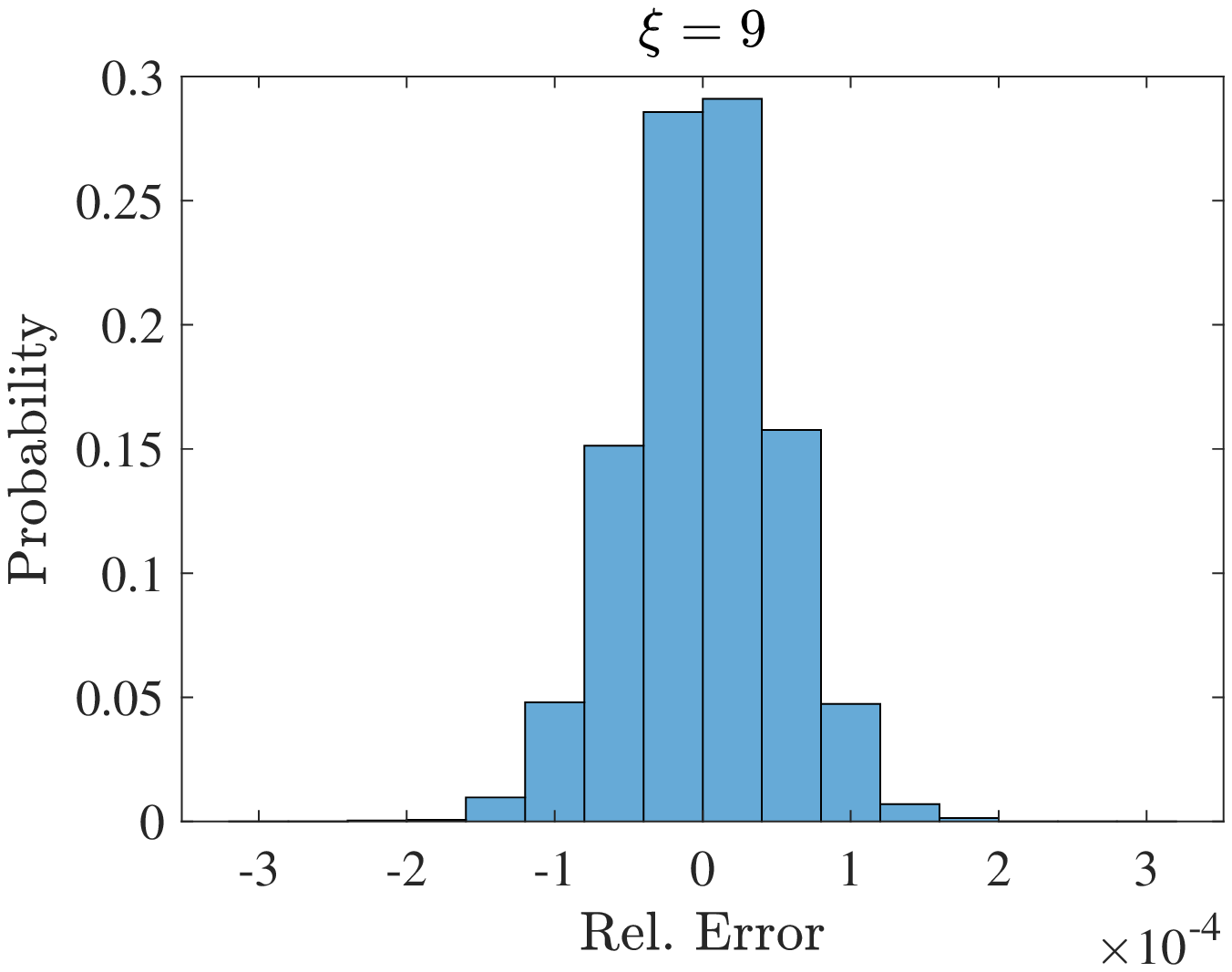}
\includegraphics[width=0.45\textwidth]{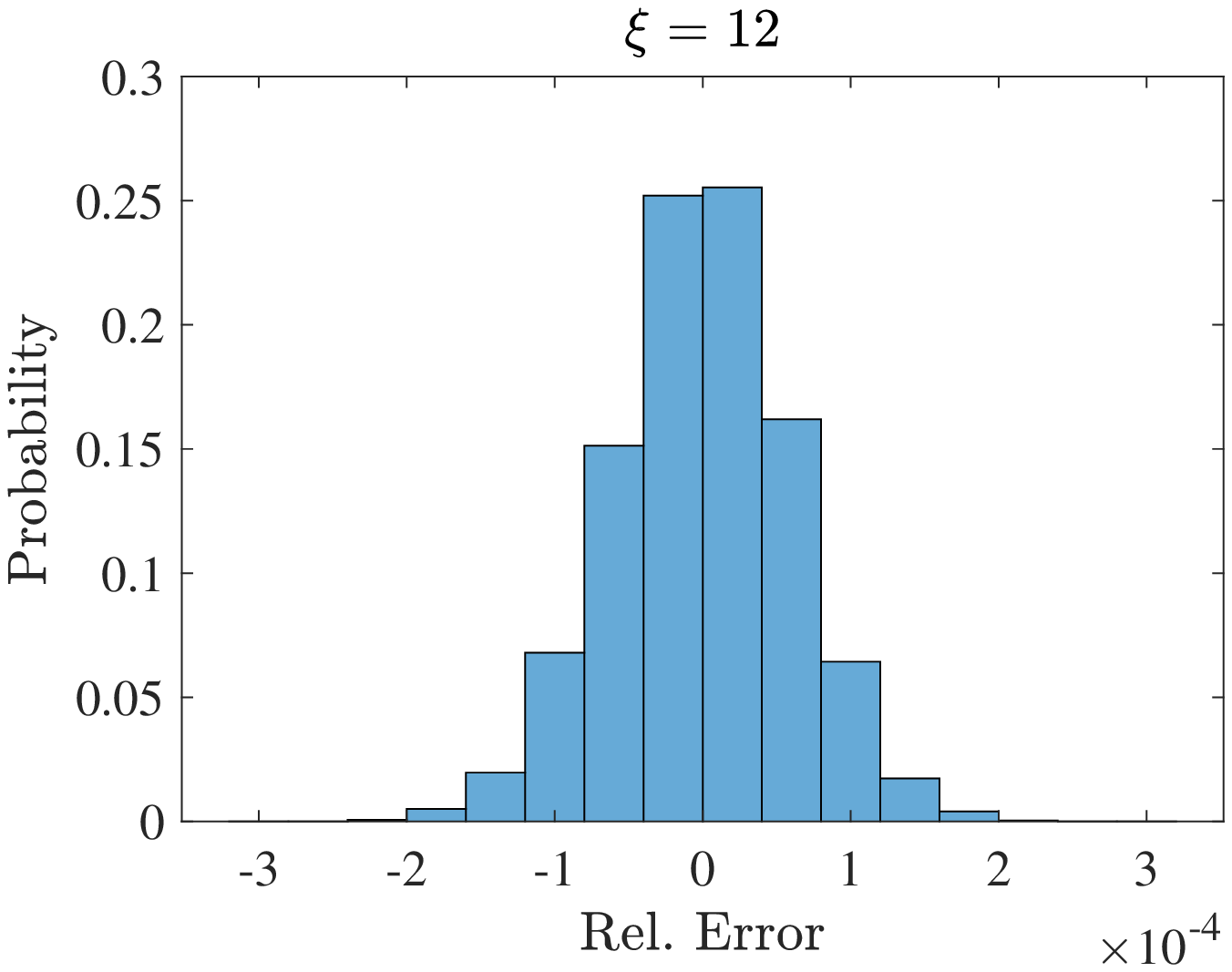}
\includegraphics[width=0.45\textwidth]{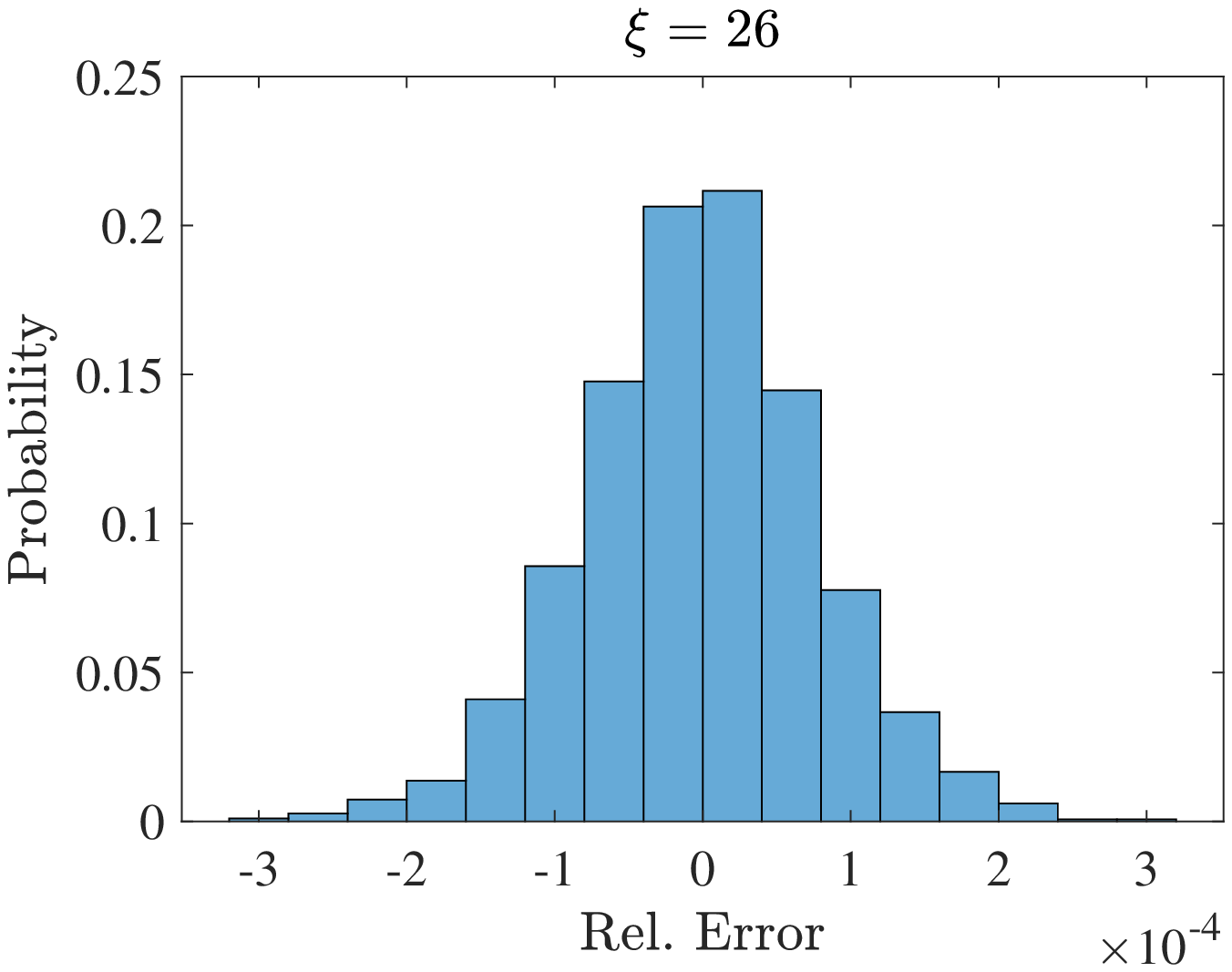}
\caption{\label{fig:ewslab} Histogram of relative errors in the components of the forces on 100 randomly positioned charges of unit strength. We compare four different values of the splitting parameter $\xi$ in the spectral Ewald method to a reference solution computed to 7 digits of accuracy using the method of Section \ref{sec:smoothdp}. The errors are normalized by the mean magnitude of the force over all charges, and the experiment is repeated 10 times. }
\end{figure}

For each Ewald parameter, we solve for the averaged potential and electric field and compute the error in each component of $\widebar{\bm{E}}(\bm{z}_i)$ (relative to the $\xi \rightarrow \infty$ reference result), normalized by the mean magnitude of $\widebar{\bm{E}}$. We do this 10 times, so that there are a total of 3000 observations $\widebar{\bm{E}}$. Figure\ \ref{fig:ewslab} shows histograms of the relative errors for each of the Ewald parameters. For all $\xi$, we obtain distributions centered approximately around 0 with 4 digits of accuracy in $\widebar{\bm{E}}$. Table\ \ref{tab:ewnoparam} shows that the standard deviation of the relative errors in each case is $\approx 5 \times 10^{-5}$. This demonstrates that the accuracy is better than the targeted $\delta = 5 \times 10^{-4}$, especially for smaller Ewald parameters when the near field dominates. 

\begin{table}
\centering
\begin{tabular}{c|c|c|c|c}
$\xi$ & Grid size & $H_E$ & $\rnf$ & Error std.\\
\hline
$\infty$ (ref.) & $128 \times 128 \times 76$ & & & \\ 
4.3 & $20 \times 20 \times 59$ & 0.5 & 0.71 & 2.7e-5\\[2 pt]
9.2 & $40 \times 40 \times 71$ & 0.25 & 0.35 & 5.1e-5\\[2 pt]
12.2 & $50 \times 50 \times 77$ & 0.20 & 0.27 & 6.0e-5\\[2 pt]
26.0 & $76 \times 76 \times 92$ & 0.13 & 0.16 & 7.9e-5
\end{tabular}
\caption{\label{tab:ewnoparam} Splitting parameters for Section\ \ref{sec:ewind}. Grid is on $[-L, L] \times [-L,L] \times [-3H_E,H+3H_E]$ with $H=0.75$ and $L=2$.  We give the standard deviation of the errors shown in Fig.\ \ref{fig:ewslab}.}
\end{table}

\subsection{Energy and force \label{sec:fge}}
Finally, we verify that the forces\ \eqref{eq:fqE} are the gradient of the energy\ \eqref{eq:Ufinal}, even in the presence of charged walls. To do this, we set $\epsilon=1$, $\epsilon_t=1/50$, $\epsilon_b=1/20$, and $H=1$ and $L=2$. We place 10 charges with unit strengths of alternating sign randomly and uniformly on $[-L,L] \times [-L,L] \times [4.5g_w,H-4.5g_w]$, where $g_w$ will be varied. We also place a Gaussian charge density on the top and bottom walls, 
\begin{equation}
    \sigma_{b/t}(x,y)=(+/-) \frac{1}{2\sqrt{4\pi^2s^4}} \exp{\left(-\frac{x^2+y^2}{2s^2}\right)}. 
\end{equation}
We use $s=0.2$ here so that the fixed charge densities can be resolved on a coarse grid, and set $\xi=6.8$ ($g_t \approx 0.07$, $H_E=0.32$, $\rnf=0.48$) for a grid of size $38 \times 38 \times 87$ points. Since the grid spacing is chosen to resolve $g_t=0.07$, and the wall charge density is a Gaussian with standard deviation $s > g_t$, we expect that our grid will  easily be able to resolve the wall density, and therefore the integrals in the second term of the energy\ \eqref{eq:Ufinal}. We have confirmed that the potential computed by the correction solve matches the analytical result for a potential of a Gaussian surface charge (not shown). 

We fix $\delta_0=10^{-4}$ and compute a rate of work by finite differences of the energy\ \eqref{eq:Ufinal}
\begin{equation}
\label{eq:w1}
W_1=-\frac{U(\bm{X}+\frac{\delta_0}{2}\delta \bm{ X})-U(\bm{X}-\frac{\delta_0}{2}\delta\bm{ X})}{\delta_0},  
\end{equation}
where $\delta \bm{ X}$ is a random displacement vector of unit length for each charge. We compare this to the rate of work done by the force on the charges\ \eqref{eq:fqE},
\begin{equation}
\label{eq:w2}
W_2=\sum_{k=1}^N \bm{F}_k \cdot \delta \bm{ X}_k = \sum_{k=1}^N q_k \widebar{\bm{E}}(\bm{z}_k) \cdot \delta \bm{ X}_k.
\end{equation}

To avoid cancellation of digits as $g_w \rightarrow 0$, we subtract the self potential, $\bar{\phi}(\bm{z}_i) \rightarrow \bar{\phi}(\bm{z}_i)-q_i\near{\widebar{G}}(r=0;g_w,\xi=0)$, so that the potentials do not become singular as $g_w \rightarrow 0$. As shown in Table \ref{tab:energslabEw}, $W_1=W_2$ to at least 3 digits even for $g_w=10^{-10}$. 

\begin{table}
\centering
\begin{tabular}{c|c|c|c|c}
$g_w$ &  $10^{-2}$ & $10^{-3}$ & $10^{-4}$ & $10^{-10}$\\
\hline
$|W_1-W_2|/W_1$ & 2.12e-4 & 3.70e-4 & 5.00e-4 & 1.40e-4
\end{tabular}
\caption{Relative differences in computed rates of work $|W_1-W_2|/W_1$ using\ \eqref{eq:w1} and \eqref{eq:w2} for $\delta_0=10^{-4}$ with $M=10$ charges.  }
\label{tab:energslabEw}
\end{table}

\section{Applications \label{sec:raul}}
In this section we apply our method to study the influence of polarization effects (image charges) on the structure of electric double layers in electrolyte solutions confined by one or two dielectric boundaries. This is a well-studied problem and it is not our intention to provide a detailed account of the topic. Rather, our main goal is to test our method by comparing to published results, and to demonstrate the importance of accounting for jumps in the dielectric permittivity. We also study the performance of the algorithm on Graphical Processing Units (GPUs) for realistic electrolyte parameters.

We use Brownian Dynamics (BD) \emph{without} hydrodynamic interactions to equilibrate an electrolyte solution in a slit channel; this requires evaluating the forces on the charges at each time step, but does not require the electrostatic energy. In a number of previous studies of electrolytes, Markov Chain Monte Carlo (MCMC) is used for equilibration. However, MCMC with local moves is expensive for electrolytes since computing the long-ranged electrostatics is a global operation and has to be done after each trial move. We have not been successful with controlling the rejection rate in MCMC methods based on global moves, including ones combining BD with MCMC like the Metropolis-adjusted Langevin Algorithm (MALA) \cite{Metropolized_MALA,MetropolizedBD}. As a compromise, we utilize here a BD method specifically designed for computing equilibrium distributions  to second-order accuracy \cite{MultistepEM_Leimkuhler}. It is important to emphasize that because we use BD, we operate in the canonical $N V T$ ensemble and not in the grand canonical ensemble. Furthermore, by construction our method requires maintaining strict electroneutrality in every configuration. Although some conventions for handling non-neutral slabs have been proposed \cite{DPPoisson_NonNeutral}, and can straightforwardly be incorporated into our handling of the $k_x=k_y=0$ mode, a unique well-defined meaning to the doubly-periodic Poisson equation cannot be given without strict electroneutrality.

We begin by summarizing our Brownian Dynamics method in Section \ref{sec:BD}. In Section \ref{sec:uncharged} we consider a monovalent binary electrolyte and study the depletion of charges near a water-air interface due to repulsion by the image charges. In Section \ref{sec:charged} we study the same effect but in a slit-channel with charged dielectric walls and with only the counterion present. Finally, in Section \ref{sec:performance} we discuss the performance of our algorithm and the optimal choice for the Ewald splitting parameter.

\subsection{\label{sec:BD}Brownian Dynamics simulations}

We consider a solution of monovalent ions with effective radius $a$ and charge $\pm e$ solvated in water. Unless stated otherwise, simulations are carried out at room temperature ($T= 298$\,K). We include a steric repulsion between ions in the form of a modified, truncated, and regularized Lennard-Jones (LJ) potential \cite{DISCOS_Periodic} with repulsive force $F_{\text{steric}}(r)=-\partial U_{\text{steric}}/\partial r$,
\begin{equation}
  \label{eq:steric}
  F_{\text{steric}}(r) =
   \begin{cases} 
     F_\text{LJ}(r=r_m) & r\leq r_m \\
     F_\text{LJ}(r) & r_m < r \leq 2^{1/p}(2a) \\
     0 & r > 2^{1/p}(2a)
   \end{cases}  
\end{equation}
where $r$ is the ion-ion distance and $F_{\text{LJ}}=-\partial U_{\text{LJ}}/\partial r$ where
\begin{equation}
\label{eq:lj}
U_{\text{LJ}} (r) = 4U_0\left(\left(\frac{2a}{r}\right)^{2p} -\left(\frac{2a}{r}\right)^{p}\right) + U_0.
\end{equation}
Here $r_m$ denotes a cutoff distance below which we molify the divergence of the classical LJ potential to avoid very large forces which lead to numerical instabilities \cite{DISCOS_Periodic}. In some simulations, we also decrease the exponent $p$ from the traditional $p=6$, in order to make the potential softer and thus increase the stable time step size. The value of the steric repulsion $U_0$ is chosen to be $\sim k_B T$ to avoid overlap but also avoid very large forces upon overlap. We add a steric repulsion with walls (dielectric boundaries) by considering a virtual image charge and computing the steric repulsion with that image charge.
 
We employ the lower precision (3 digits of accuracy) set of parameters for our algorithm (see last column in Table \ref{tab:ewconstraints}); we have confirmed that the forces on the charges did not change to at least 3 digits upon changing the Ewald splitting paramater $\xi$. The algorithm is implemented on GPUs using the Universally Adaptable Multiscale Molecular Dynamics (UAMMD) \cite{uammd} framework, and we use single precision in all calculations, except that the correction potential\ \eqref{eq:inphisc} is computed in double precision to avoid overflow and underflow. We have confirmed that using single precision in the rest of the calculations does not lead to a degradation of the overall accuracy of the method. Our GPU code is freely available on the UAMMD github \url{https://github.com/RaulPPelaez/UAMMD}, and the tests given in this section can be reproduced using the codes and input files available at  \url{https://github.com/stochasticHydroTools/DPPoissonTests/}.

We solve the stochastic equations of Brownian Dynamics without hydrodynamic interactions,
\begin{equation}
d\V{z} = \mu \left( \V{F}\left(\V{z}\right) + \V{F}_{\text{steric}}\left(\V{z}\right) \right) dt + \sqrt{2 \kT \mu}\; d\V{B}
\end{equation} 
where the electrostatic forces $\V{F}$ are defined in \eqref{eq:fqE}, $\mu$ is a scalar mobility, and $\V{B}(t)$ is a collection of independent Brownian motions (Wiener processes). These equations are integrated in time using a non-Markov variation of the classical Euler-Maruyama method proposed in \cite{leimkuhler} and studied theoretically in \cite{MultistepEM_Leimkuhler}. The scheme updates the ion positions from time step $n$ to $n+1$ using
\begin{equation}
\label{eq:leimkuhler}
\V{z}^{n+1} =\V{z}^{n} + \mu \left( \V{F}_{\text{steric}}\left(\V{z}^n\right) + \V{F}\left(\V{z}^n\right) \right) \D{t} + \sqrt{\frac{\kT \mu \D{t}}{2}} \left(\V{W}^{n} + \V{W}^{n+1}\right),  
\end{equation}
where $\V{W}^n$ is a collection of independent Gaussian random variables of mean zero and variance one. The scheme \eqref{eq:leimkuhler} has first order weak accuracy at finite times but transitions to second-order weak accuracy exponentially fast in time \cite{MultistepEM_Leimkuhler}. In particular, it approximates expectation values with respect to the Gibbs-Boltzmann equilibrium distribution $\sim \exp\left(-\left(U+U_{\text{steric}}\right)/\left(\kT\right)\right)$ to second-order accuracy in the time step size $\D{t}$. We have found this scheme to be the most accurate among several we tested, and also found it to be just as robust for larger time step sizes as the less-accurate Euler-Maruyama scheme. We report the values of $\D{t}$ and time in units of the diffusive time $\tau_0= a^2/\left(\kT \mu\right)$, so that the precise value of $\mu$ is irrelevant.
Charges are checked every step to be inside the allowed domain (between $z=0$ to $z=H$). If a charge leaves the domain after a given update the step is discarded and retried until it results in a valid configuration. After a maximum allowed number of retries an unrecoverable configuration is declared to avoid stagnation. We choose $\D{t}$ to be small enough such that steps are rejected very rarely and unrecoverable configurations do not occur over long periods of time. We find that we can use a larger time step size if we restrict the maximum allowed displacement per time step for a particle to a multiple of the ion radius; we return to this point in the Conclusions.

\begin{figure}
\includegraphics[width=\linewidth]{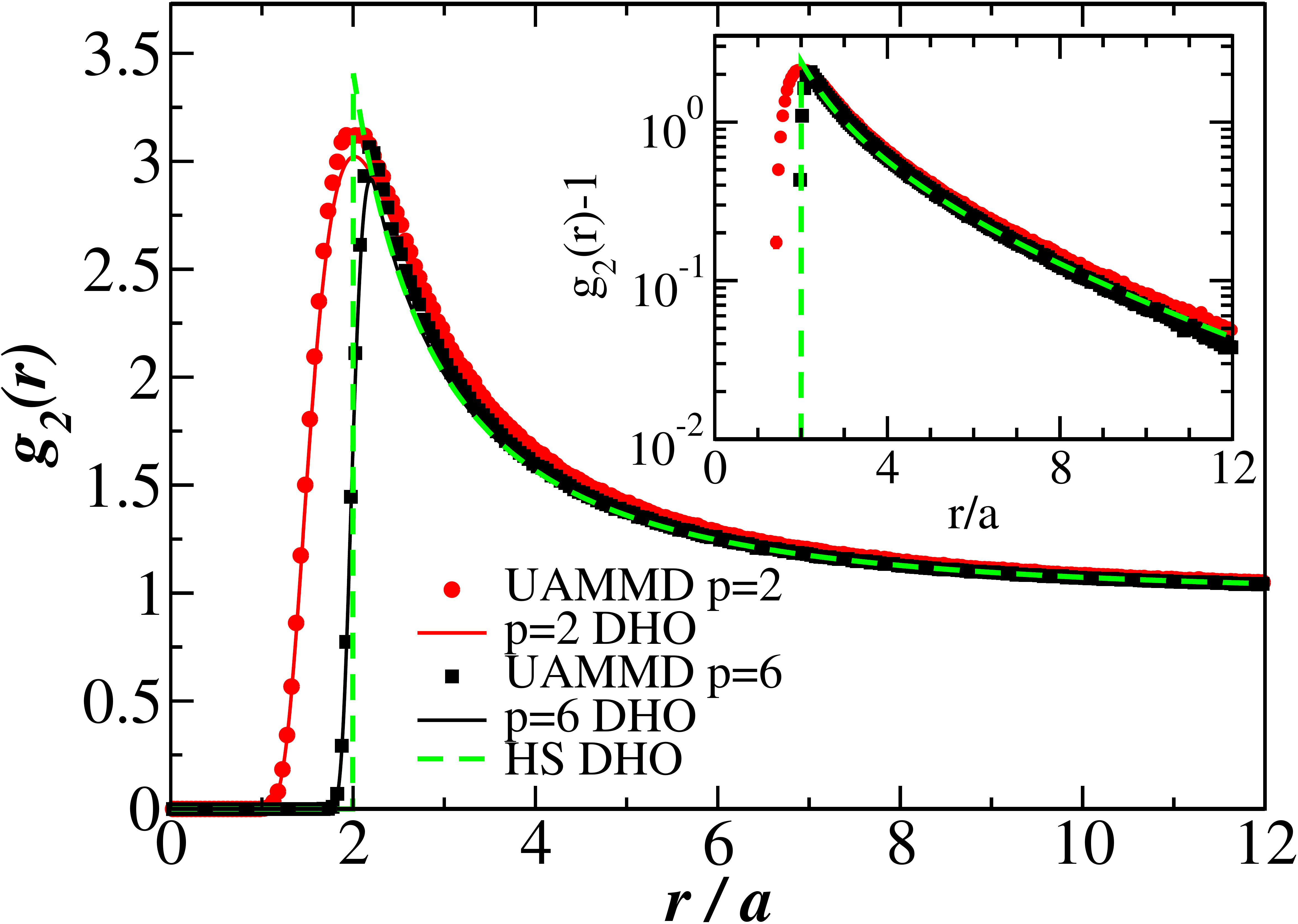}
\caption{\label{fig:gdrp6}Pair correlation function (symbols) of counter ions in a dilute binary monovalent electrolyte of molarity 0.05M in a triply periodic cubic domain ($N=8192$ ions), for two different systems ($p=6$ or $p=2$ in the LJ potential \eqref{eq:lj}). The parameters for the potential ($p$, $r_m$, $U_0$, and $g_w$) and time step size $\D{t}$ are given in Table \ref{tab:parameters}. Solid lines are the DHO theory \eqref{eq:gdrdho}, and the dashed line shows the theory for a hard-sphere steric repulsion with the rest of the parameters as for the system with $p=6$. The inset shows the exponential tail $\sim\exp(-r/\lambda)/r$ in $g_2(r)-1$.}
\end{figure}

It is straightforward to adapt our method and implementation from a doubly-periodic to a triply-periodic domain. In Fig.\ \ref{fig:gdrp6} we show results for the pair correlation function $g_2(r)$ of counter ions in a dilute binary monovalent electrolyte of molarity 0.05M in a triply periodic domain. The parameters for the simulations are given in Table \ref{tab:parameters} and in the caption of the figure. In particular, we study two different systems in the remainder of this section. The first system (second column of table) is meant to correspond to the hard-sphere system studied theoretically and using Monte Carlo simulations in \cite{mcmc1980}. To mimic hard spheres, we set the exponent to $p=6$ in the LJ potential \eqref{eq:lj} and choose $U_0$ to be sufficiently large to make $g_2(r)$ almost zero for $r<2 a$. The second system (third column of table) we study uses a softer repulsion with exponent $p=2$ to avoid stiffness. In this second system, however, we wish to mimic point charges closely so we set $g_w \ll a$, which leads to stiffness of the electrostatic potential. When neither strict steric exclusion or point-charge electrostatics is required, one can use a small $p=2$, large $r_m \sim a$, and large $g_w \sim a$, to make the combined interaction potential non-stiff and allow for a larger $\D{t}$; we return to this point in the Conclusions.

\begin{table}
	\centering
	\begin{tabular}{ c | c | c}
		Figure & \ref{fig:mcmc}& \ref{fig:charged_wall}\\
		\hline
		$a$& $2.125\angstrom$ & $ 2\angstrom $\\
		$L_{xy} $&$ 185.8a$ & $800a$\\
		$H$&$50a$&$100a$\\
		$g_w$&$0.25a$&$0.02a$\\
		$U_0$&$0.7005\kT$&$0.2233\kT$\\
		$r_m$&$1.5a$&$a$\\
		$p$&$6$&$2$\\
	\end{tabular}
	\caption{\label{tab:parameters}Parameters used in the BD simulations corresponding to Figs.\ \ref{fig:mcmc} and\ \ref{fig:charged_wall}. The time step size is $\Delta t = 2 \cdot 10^{-3} \tau_0$.}
\end{table}

For comparison to the numerical results, in Fig. \ref{fig:gdrp6} we also show theoretical predictions for the pair correlation function based on Debye-Huckel-Onsager theory, modified to account for the fact that our charges are not point charges,
\begin{equation}
\label{eq:gdrdho}
g_\text{DHO}(r) = \exp\left(-\ddfrac{U_{\text{steric}}(r) + U(r)\exp(-r/\lambda)}{\kT}\right)
\end{equation}
where $U(r)$ is the electrostatic potential between two counter ions with Gaussian charge density (see first term in \eqref{eq:GnearK}),
\begin{equation}
\label{eq:gaussian}
U(r) = -\ddfrac{e^2}{4\pi\epsilon r}\;\text{erf}\left(\ddfrac{r}{2 g_w}\right).
\end{equation}
Here the Debye length is
\begin{equation}
\label{eq:debye}
\lambda = \left(\ddfrac{\epsilon \kT}{e^2N_a(2 M)}\right)^{\frac{1}{2}},
\end{equation}
where $N_a$ is the Avogadro number and $M$ is the salt molar concentration (molarity) in  $\text{mol}/m^3$. The numerical results in Fig. \ref{fig:gdrp6} are in good agreement with the theory especially in the Debye (exponential) tail (see inset). For the system with $p=6$ we also show the theoretical $g_2(r)$ if the steric repulsion were a hard-sphere potential (strictly no overlaps).

\subsection{\label{sec:uncharged}Uncharged wall}

\begin{figure}
  \centering
  \includegraphics[width=\linewidth]{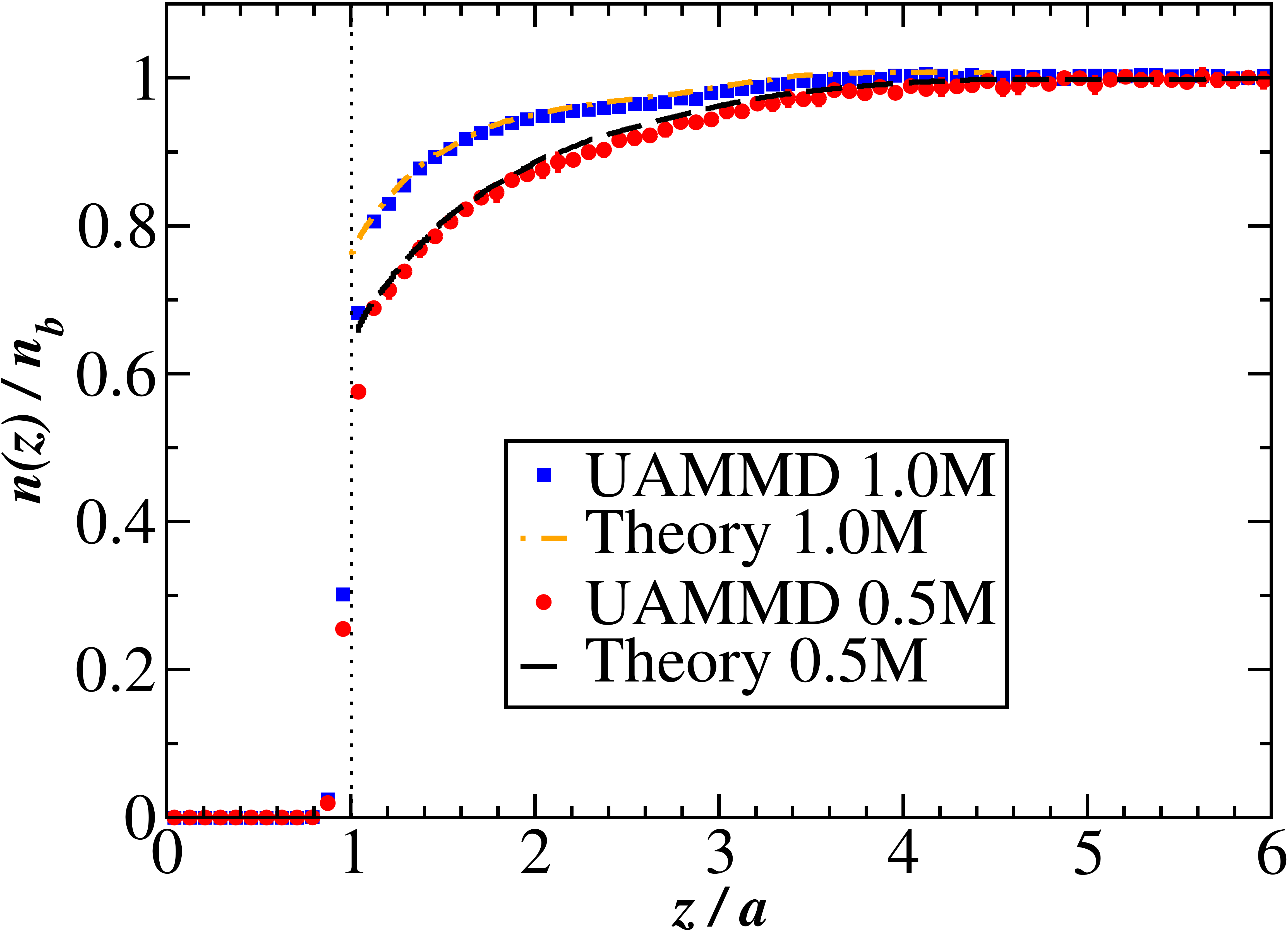}
\caption{\label{fig:mcmc}Number density of monovalent ions $n(z)$ at a distance $z$ from a dielectric boundary, normalized by the ``bulk'' density far away from the wall. The results from our code UAMMD are compared with the results in figures 3 and 4 in \cite{mcmc1980}. The bulk ionic density is $n_b= 0.495$ mol/L $= 2.8\cdot 10^{-3}a^{-3}$ (black line) or $n_b = 1.00$ mol/L $= 5.8\cdot 10^{-3}a^{-3}$ (orange line); other parameters are listed in the second column in Table \ref{tab:parameters}. The total number of ions is $N=9398$ (red dots) and $N=19142$ (blue dots). Simulations were initialized by randomly placing ions in the region $a<z<H-a$ and equilibrated for 20 diffusive times $\tau_0$, and then data was collected every diffusion time for the subsequent 200 diffusive times. Error bars (comparable to or smaller than the symbol size) come from averaging over 32 independent simulations.}
\end{figure}

In this subsection we reproduce results from the article \cite{mcmc1980} using the doubly-periodic Poisson solver in UAMMD.
We simulate an electroneutral monovalent solution inside a slab with uncharged walls and different permittivities above and below it; the parameters are listed in the second column in table \ref{tab:parameters}. We compare our results and those obtained in \cite{mcmc1980} in Fig.\ \ref{fig:mcmc}. Since the walls are uncharged there is a symmetry between the anions and cations, and therefore we show the ion number density $n(z)$ averaged among the two types of ions.

There are a number of recent theories that account for the polarization effects next to a dielectric boundary, but all of them involve nontrivial computations --- for recent examples see \cite{ModifiedPNP_Polarization} without and \cite{ModifiedPNP_Polarization_HardSphere} with steric repulsion. The authors of \cite{mcmc1980} developed several theoretical approaches to compute the equilibrium density of ions next to a planar interface (wall) with hard-sphere steric repulsion in addition to electrostatics (i.e., the so-called primitive model of electrolytes), with the ``BBGKY+EN'' theory being in best agreement for the molarities we study here. We have extracted the theoretical curves from the figures in \cite{mcmc1980} and show those for comparison in Fig.\ \ref{fig:mcmc}. 

The MCMC simulations and theory in \cite{mcmc1980} used the grand canonical ensemble and considered a system that has a bottom wall but is unbounded in the $z$ direction, i.e., a system that is in contact with an infinite reservoir as $z/a\rightarrow\infty$ with bulk molarity $M_{\text{bulk}}$. Since we use the canonical ensemble, we put a top wall at $z=H=50 a$ but set $\epsilon_t=\epsilon$ in order to mimic an ``open'' reservoir boundary. We numerically integrated the theoretical curves in \cite{mcmc1980} to obtain a total number of ions that would give the same bulk molarity far away from the wall, assuming that the density was constant for $z \gg a$ (which is in agreement with simulation results).

The results in Fig. \ref{fig:mcmc} show excellent agreement between our numerical results and the theoretical predictions from  \cite{mcmc1980}, both for bulk molarity $M_{\text{bulk}}\approx 0.5$M and $M_{\text{bulk}}\approx 1$M, validating our method and implementation. We see that there is a depletion of ions next to the interface because of the repulsion with the image charges, with some modest layering due to steric repulsion for $M_{\text{bulk}}=1$M predicted by the theory, and confirmed by the simulations. In this example, the depletion layer extends only over a distance of a couple of ion diameters away from the wall.

\subsection{\label{sec:charged}Charged wall test}

In this subsection we study a slit channel with charged walls with only the counterions present in the interior of the channel. Specifically, we place $N$ ions of charge $+e$ in the channel and keep the system electroneutral by setting the surface charge densities at the top and bottom walls to $\sigma_b = \sigma_t = \sigma = -(N e)/\left(2L_{xy}^2\right)$; the parameters are listed in the third column in Table \ref{tab:parameters}.  For this special case it is possible to analytically solve the Poisson-Nernst-Planck (PNP) equations for the equilibrium number density of ions across the channel (see Section 4.1 in \cite{ElectrostaticsMembranes}).
If we denote the sterically-accessible thickness of the channel as $d=H-2a$ and solve the PNP equations with charged walls inside a channel of width $d$, we obtain the ion density profile for $a<z<H-a$ (zero outside of this interval),
\begin{equation}
\label{eq:pnp}
n(z) = \ddfrac{n_m}{\cos\left(K (z-H/2)/d\right)^2},
\end{equation}
where the dimensionless constant $K<\pi$ depends on the charge density,
\begin{equation}
\sigma = -\frac{e}{2}\int_{z=a}^{H-a} n(z)dz = -\ddfrac{2 K \tan\left(K/2\right) \epsilon \kT}{d e},
\end{equation}
and $n_m$ is the number density of ions in the middle of the channel,
\begin{equation}
\label{eq:nm}
n_m =\frac{2K^2\epsilon \kT}{d^2 e^2}.
\end{equation}

\begin{figure}
  \centering
  \includegraphics[width=0.49\linewidth]{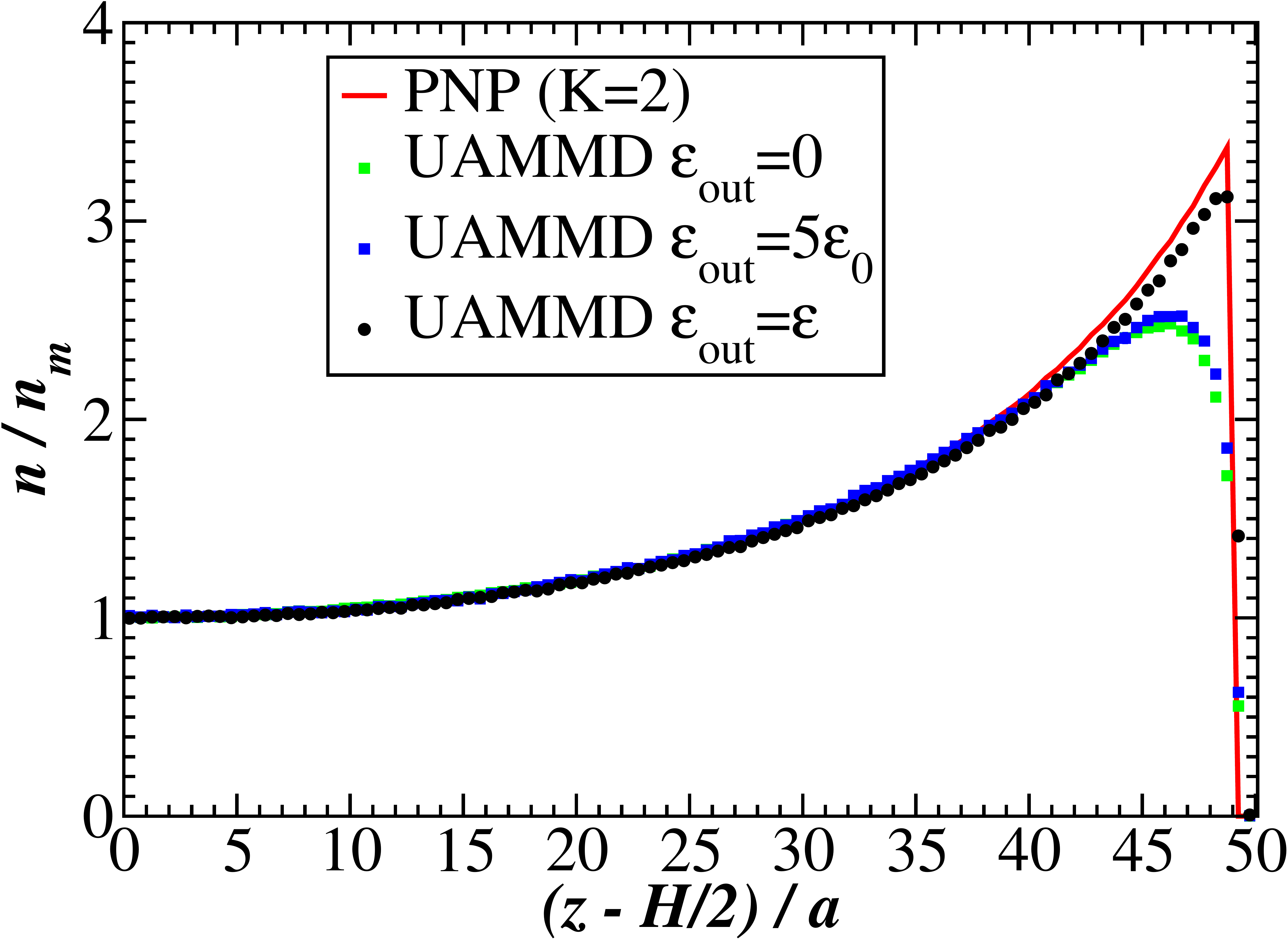}
  \label{fig:charged_denser}\includegraphics[width=0.49\linewidth]{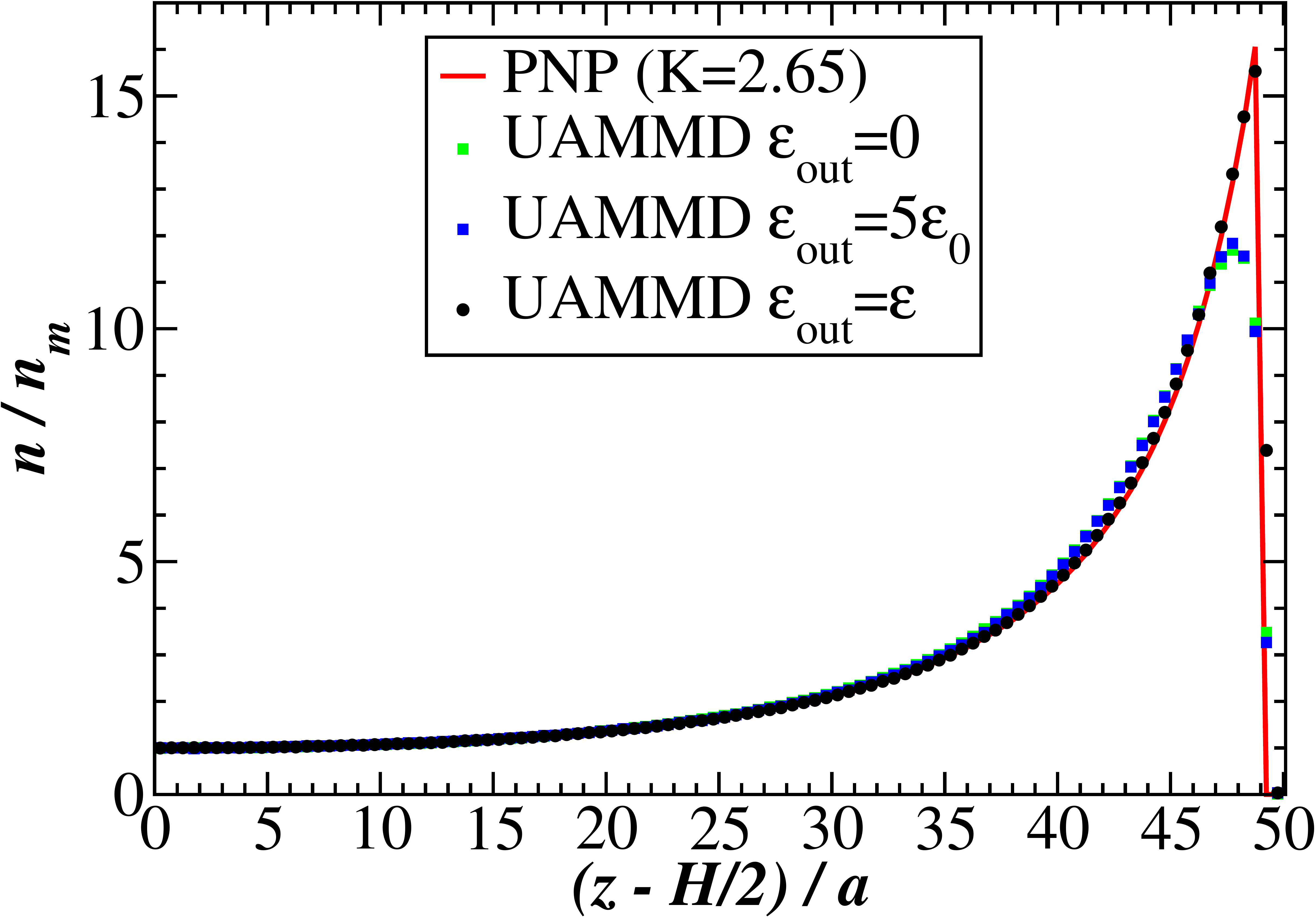}
\caption{\label{fig:charged_wall}Number density of cations $n(z)$ inside a slab channel with oppositely charged walls, at two different surface change densities (lower in the left panel). The density is normalized by the density in the middle of the channel predicted by the PNP equations (see \eqref{eq:nm}), $n_m = 3.26\cdot 10^{-5} a^{-3}$ (left, $N=6140$ ions) and $n_m = 1.85\cdot 10^{-5} a^{-3}$ (right, $N=1739$ ions); other parameters are listed in the third column in Table \ref{tab:parameters}. Numerical results obtained using the UAMMD doubly-periodic code are shown for three different permittivities outside of the channel (glass, water, and an unphysical zero permittivity medium). To accelerate equilibration, initial configurations were generated by randomly and uniformly placing ions and using rejection to make the density vary with $z$ according to \eqref{eq:nm}. We equilibrated for 20 diffusive times $\tau_0$, and then data was collected every diffusion time for the subsequent 200 diffusive times. Results are averaged over $64$ independent runs; the error bars are smaller than the symbol size.}
\end{figure}

In Fig.\ \ref{fig:charged_wall} we compare numerical results for $n(z)$ with the prediction of the PNP equations \eqref{eq:pnp} at two different densities. We see an accumulation of ions near the boundaries because of the surface charge, but the maximum density depends on the dielectric contrast. We consider three different materials on the outside of the channel: glass (permittivity $\approx 5$ times that of vacuum), water (no dielectric jump), or an unphysical medium with zero permittivity.

We find good agreement with the PNP prediction only if there is no dielectric jump and $\epsilon_t=\epsilon_b=\epsilon_{\text{out}}=\epsilon$. This is not unexpected since the PNP equation does not take into account polarization effects. Namely, for an electroneutral slit (or cylindrical) channel, Gauss's law dictates that the electric field just outside of the channel must be zero; the electric fields created by the ions on the outside of the channel are evanescent and vanish in the mean field approximation. Only if charge electroneutrality is broken does the dielectric constant on the outside of the channel begin to matter in the PNP equations \cite{PNP_NonElectroneutral}. While polarization effects can be included in modified PNP equations that take into account charge correlations \cite{ModifiedPNP_Polarization,ModifiedPNP_Polarization_HardSphere}, the resulting nonlocal equations are no longer simple to solve even numerically. It is important to note that the mismatch between the PNP theory and the numerical results seen in Fig.\ \ref{fig:charged_wall} extends over a relatively large distance from the wall, on the order of tens of particle radii, emphasizing the need to account for dielectric jumps in theoretical and numerical descriptions of confined electrolytes.

The results in Fig.\ \ref{fig:charged_wall} show that there is only a small difference between the results obtained for glass and a fictitious material with $\epsilon_{\text{out}}=0$, for which there is no dielectric displacement outside of the channel. This is good news for approaches such as immersed boundary methods \cite{DISCOS_Periodic} that only resolve the electric fields inside the channel, and impose (inhomogeneous) Neumann conditions on the electric field, rather than the correct jump boundary conditions for the dielectric displacement.

\subsection{\label{sec:performance}Computational performance}

In this section we study the computational performance of our GPU implementation of the doubly-periodic (DP) method, as a function of the Ewald splitting parameter $\xi$. We randomly and uniformly place $2\cdot10^4$ charges inside a slab of extent $H=50$ and length $L_{xy}=185$, and set $g_w=a/4$. We use single precision (except for the correction solve to avoid overflow) and time each part of the algorithm using profiling tools; this requires that all kernel launches are synchronous, not allowing for overlapping CPU/GPU work and concurrent kernels (which improve performance only slightly). 

\begin{figure}
  \centering
  \includegraphics[width=0.49\linewidth]{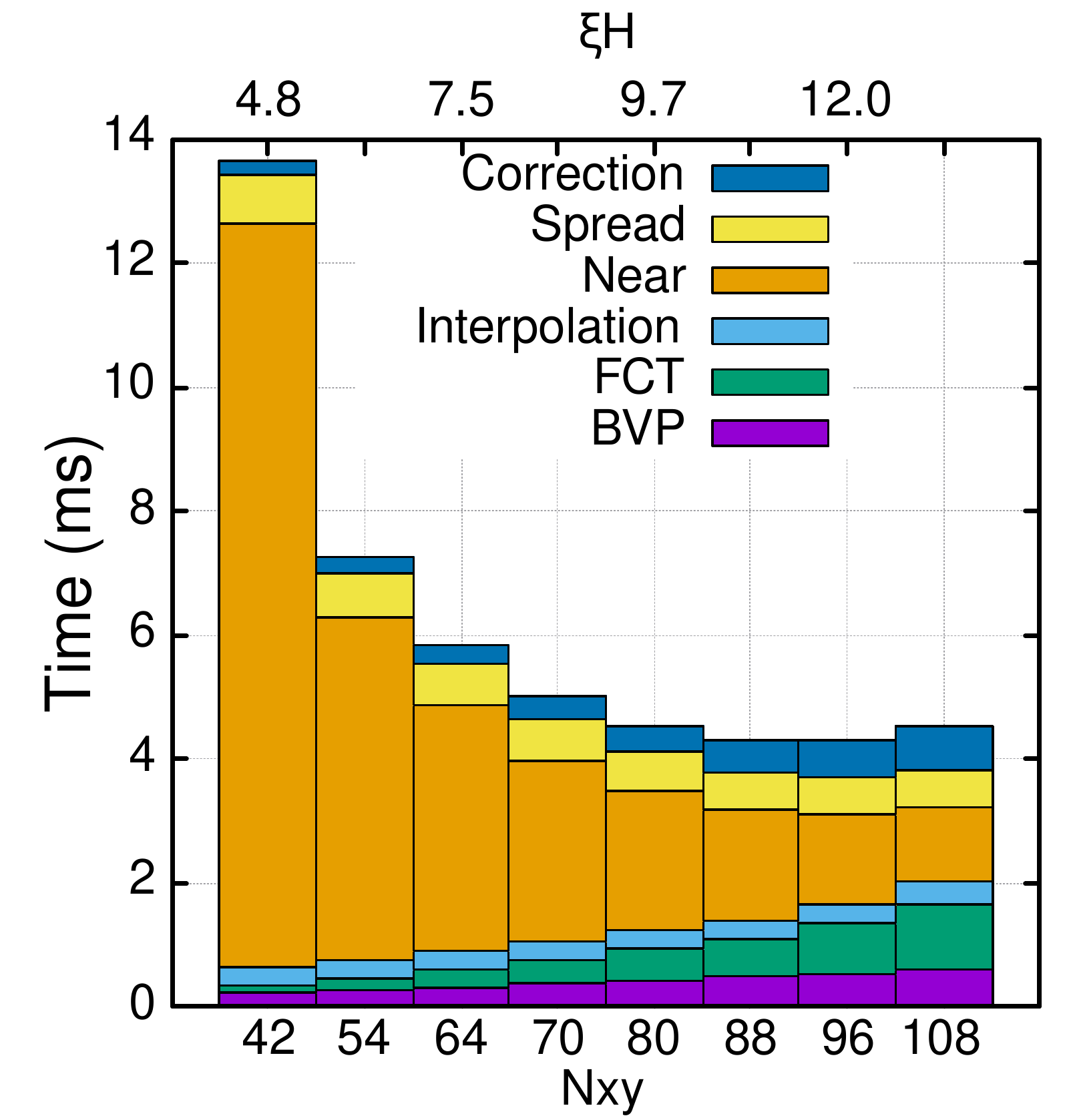}
  \includegraphics[width=0.49\linewidth]{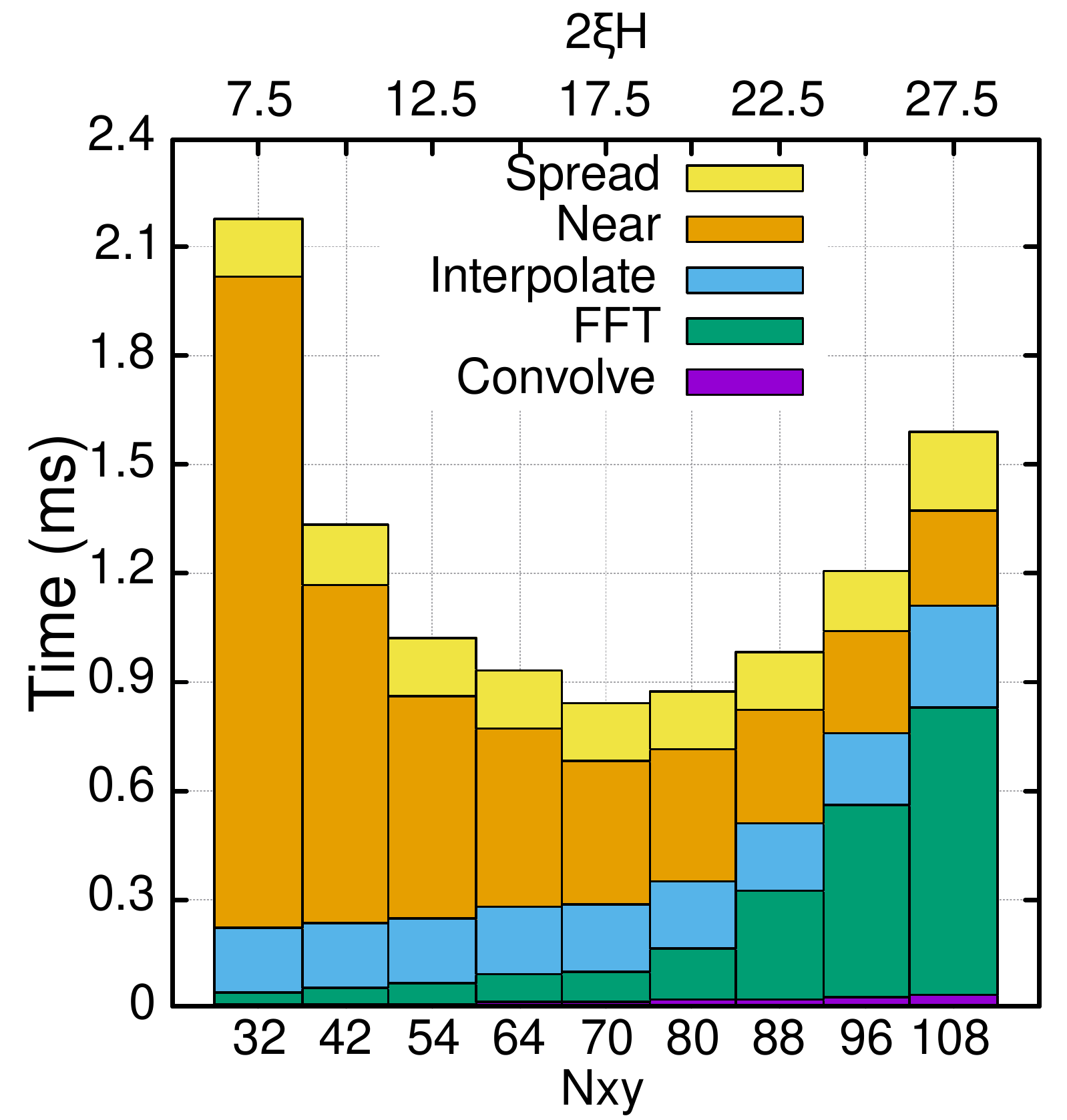}
\caption{\label{fig:perf}Execution time for different components of the algorithm for $2\cdot10^4$ charges as a function of the Ewald splitting parameter $\xi$ (top axis), which controls the grid size in the $xy$ directions (bottom axis). The left panel is for a doubly-periodic (DP) domain, while the right one is for a triply-periodic (TP) domain. Timings are collected using the UAMMD \cite{uammd} code compiled and run on a RTX2080Ti GPU using NVIDIA's CUDA 11.0.}
\end{figure}

The total wall clock time to compute the force on and potential at each charge, and the breakup among different components of the algorithm, are shown in the left panel of Fig.\ \ref{fig:perf}. We isolate the time due to the following components of the algorithm: Computing the correction field to account for the missing images, spreading the charges to the grid, computing the near-field interactions, interpolating the field back onto the charges, computing the Fourier-Chebyshev Transform (FCT) using 3D FFTs, and solving the boundary value problems (BVPs) in the doubly-periodic Poisson solver.

For comparison, in the right panel of Fig.\ \ref{fig:perf} we also give timing results for a triply periodic (TP) domain (also available in UAMMD \cite{uammd}), for which the Poisson equation can be solved entirely in Fourier space using 3D FFTs. We use the same particle configurations as for DP domains in the TP domain, but we set the length of the domain in the $z$ direction to $L_z=2H$, which makes the sizes of the grids used in the 3D FFTs similar. Unlike the DP solver, the TP solver does not require any images, a correction solve, or BVP solver, and, importantly, uses a uniform grid for which spreading and interpolation are faster in the $z$ direction.\footnote{In particular, fast Gaussian gridding \cite{NUFFT} can be used in all directions, although we do not find this to lead to substantial improvement on the GPU. On the GPU, the main savings comes from the fact that the number of grid points to spread to or interpolate from in the $z$ direction is constant (10 grid cells in each direction) in the TP case.}

At the optimum split (grid size $N_{xy} = 88$ for DP and $70$ for TP), the time it takes to complete the calculation is $0.84$ ms for TP and $4.3$ ms for DP. While in this example we find that DP is a factor of $5$ slower than TP, the total execution time is still remarkably small given that no supercomputer hardware is required; the GPU used is widely available and inexpensive. Note that other parts of the algorithm such as Brownian Dynamics updates and computing the steric interactions take only a negligible fraction of the time it takes to compute the long-ranged interactions.

\section{Conclusions \label{sec:conclusions}}
We have developed a spectrally-accurate fast method for computing electrostatic energy and forces for a collection of charges in doubly-periodic slabs with jumps in the dielectric permittivity at the slab boundaries. To do this, we used a modification of the Spectral Ewald (SE) method to smear point-like charges into Gaussian clouds. Unlike existing methods based on Fourier transforms in all directions, our method uses a novel Fourier-Chebyshev solver to solve the smoothed (far-field) Poisson equation, which is restricted to a finite domain using the Dirichlet-to-Neumann map.  The grid size for this solver is controlled by the Ewald splitting parameter and can be optimized to balance near-field and far-field costs for optimal performance. We used an image construction to handle Ewald charge clouds that overlap the dielectric boundaries, and handled the remaining \emph{smooth} mismatch in the boundary conditions (BCs) at the dielectric interfaces using the analytic solution of the Laplace equation with inhomogeneous boundary conditions. Combining the far-field and correction steps yields a method that requires the same components as the Spectral Ewald method for triply-periodic domains: spreading and interpolation from the charges to the grid using Gaussian kernels (but note that the Chebyshev grid is \emph{not} uniform in the $z$ direction), and a forward and inverse \emph{three-dimensional} FFT. All components can be parallelized on GPUs, and our public-domain implementation takes only $\sim5$ms per time step of Brownian Dynamics (BD) for 20K charges in a model electrolyte solution.

We used the GPU BD code to study the equilibrium structure of the (Debye) double layer next to an uncharged and a charged dielectric interface between water (the solvent) and a material with a low dielectric permittivity (like glass or air). We found that the interaction with the image charges causes a substantial depletion of charges next to the interface over a layer that can extend many ion diameters. This effect is not predicted by mean-field PNP theories but has to be accounted for whenever there is a substantial dielectric jump.
In many previous studies, a number of uncontrolled approximations are made that can only be justified by comparing to an algorithm that does not make those approximations. For example, it is well-known that in binary electrolytes electrostatic interactions are screened by the Debye counterion cloud. This allows some authors \cite{mcmc1980} to neglect periodic effects in the $xy$ direction and thus compute electrostatic interactions by summing over only the nearest image when computing electrostatic energy. Our results in Fig. \ref{fig:mcmc} appear to be in good agreement with the results given in \cite{mcmc1980}, suggesting this approximation was justified for the system studied. However, we also study here a case where there are no counterions, and there is no screening, so that it is crucial to properly implement periodicity in the $xy$ directions. 

Another important approximation that greatly simplifies the problem is to assume that the dielectric constant is zero outside of the slit channel. While this is unphysical (the smallest possible permittivity is that of vacuum), the dielectric constant of most materials such as glass or lipid membranes is much lower than that of water, which is the typical solvent. In the unphysical limit of zero dielectric constant in the exterior of the channel, instead of jump conditions on the electric displacement we get a Neumann condition on the electrostatic potential on the interior of the channel. This means that it is no longer necessary to compute the electrostatic potential outside of the channel or worry about infinitely many images, similar to the case of metallic boundaries. Our results in Fig. \ref{fig:charged_wall} suggest that for a water-glass interface the approximation of zero permittivity outside the channel is quite accurate.

It is straightforward to incorporate metallic slab boundaries (electrodes) in our methods instead of dielectric jumps; only minor modifications to the correction solve are required to set the potential at the top and/or bottom electrodes to a specified value. Since the case of two dielectric jumps is the hardest and requires the full power of our method, we focused on this case in our tests.

Our experience with BD for electrolyte solutions suggests that an important problem to overcome in future work is the small time step size required to stably integrate a system with stiff steric repulsion. To mimic hard-sphere repulsion in the presence of walls, we required $\D{t}\approx 2\cdot10^{-3}\tau_0$, where $\tau_0$ is the typical time it takes an ion to diffuse a distance of its radius. Even with the algorithmic improvements we developed here and the computational power of GPUs, this $\D{t}$ is too small to reach time scales relevant to dynamics of electrolytes, including relaxation and charging dynamics of double layers, and electrohydrodynamic phenomena. As in recent work based on the Discrete Ion Stochastic Continuum Overdamped Solvent (DISCOS) method \cite{DISCOS_Periodic}, here we introduced several mollifications that helped increase $\D{t}$. The first one was to mollify the charges by making them Gaussian clouds of finite width instead of point charges, which helps avoid the $r^{-2}$ divergence of electrostatic forces for pairs of ions. We also mollified the traditional Lenard-Jones steric repulsion by capping its divergence at overlap and softening the potential by lowering the exponent to $p<6$. These two changes are related to each other and have to be done in tandem because the steric repulsion (mimicking Pauli exclusion) has to be strong enough to prevent overlap of counter ions. 

The primary source of instability in temporal integration appears to be large steric or electrostatic forces that occur upon ion overlap. These large forces occasionally lead to displacements that are several ionic radii large, which is particularly problematic for ions near walls, as the displacements can lead to ions leaving the domain. We have had some success increasing the time step size by limiting the largest possible displacement of an ion during a time step to a fraction of the ionic radius. Even with this ad hoc change the largest stable and accurate time step size we achieved was $\D{t}\approx 10^{-2}\tau_0$, which is still quite small even though it is much larger than $\D{t}$ in MD. We hope that future mathematical study of temporal integrators for overdamped BD with stiff steric repulsion will lead to improvements and allow us to reach a desirable $\D{t}\approx 10^{-1}\tau_0$. Another avenue worthy of exploration is reusing information between time steps, especially for the far-field Poisson solve. Since the far-field potential and fields are smooth, it is perhaps possible to not repeat all steps of the solve at each time step. Our preliminary investigations found some promise in this direction for dilute electrolytes.

Another important direction is to generalize our method to Stokes flow so that hydrodynamic interactions can be accounted for in Brownian Dynamics. Our doubly-periodic Fourier-Chebyshev solver can straightforwardly be generalized to the Stokes instead of the Poisson equation, as we will present in future publications. The main challenge is Ewald splitting in the presence of no-slip boundaries (bottom wall only or top and bottom walls); for triply periodic systems one can use the Positively Split Ewald method \cite{SpectralRPY} or related SE methods \cite{SpectralEwald_Stokes} that rely heavily on Fourier transforms in all directions. Some progress on real-space based Ewald splitting with boundaries has been made in \cite{BrownianDynamics_OrderN}; however, because an image construction was not used to handle the boundaries, the near field does not satisfy the BCs on the wall (as it did in our method for the Poisson equation). While the mismatch in BCs can in principle be fixed with a correction solve \cite{BrownianDynamics_OrderN}, the grid required for an accurate correction would be much finer than the grid used for the far-field solver, negating the advantages of Ewald splitting. Furthermore, the method in \cite{BrownianDynamics_OrderN} cannot handle a single bottom wall as does our approach based on the Dirichlet to Neumann map. 

Ewald methods for Stokes flow based on image constructions have been developed for a single bottom wall using Fourier transforms in all directions \cite{SpectralEwald_Wall} or FMMs \cite{FMM_wall,STKFMM}. However, the image construction for a no-slip wall for Stokes flow involves several types of image singularities \cite{OseenBlake_FMM}, and this leads to substantial complexity and inefficiency compared to the approach we developed here for the Poisson equation. It should be mentioned that recent investigations using the DISCOS method demonstrate that hydrodynamic interactions make important contributions to transport in electrolytes \cite{DISCOS_Periodic}. That said, these recent studies also demonstrate that hydrodynamic interactions can be coarse grained at scales smaller than the typical ion-ion distance and replaced by standard non-hydrodynamic or ``dry'' diffusion (as we used in this work), which suggests that Ewald splitting may not be necessary for electrolyte solutions. Nevertheless, it remains a challenge for the future to adapt the method developed here to Stokes flow.

\begin{acknowledgments}
We thank Zecheng Gan for helpful discussions regarding electrostatic energy in the presence of surface charges. Ondrej Maxian is supported by the National Science Foundation (NSF) via GRFP/DGE-1342536. This work was also supported by the NSF under award DMS-2011544 and through a Research and Training Group in Modeling and Simulation under award RTG/DMS-1646339. Ra\'ul P. Pel\'aez acknowledges funding from Spanish government MINECO project FIS2017-86007-C3-1, and thanks Prof. Rafael Delgado-Buscalioni for his support and additional funding.
\end{acknowledgments}

\textbf{Data availability:} All of the codes and input files to reproduce our results are freely available at
\url{https://github.com/stochasticHydroTools/DPPoissonTests}.

\appendix
\section{Boundary value solver \label{sec:bvps}}
Our boundary value problem (BVP) solver for\ \eqref{eq:BVP} is based on the specgtral integration method of \cite{greengard1991spectral}. Without loss of generality, it is most convenient when reviewing this formulation to assume the domain is $z \in [-1,1]$. In this case, the BVP we need to solve is generally of the form
\begin{gather}
y''(z)-k^2y(z)=f(z),\\[2 pt]
y'(1)+ ky(1)=\alpha, \qquad y'(-1)- ky(-1)=\beta.
\end{gather}
An analytical solution can be derived for this BVP, but it requires numerically computing an integral with integrand related to $e^{kz}$. Since this calculation must be done using very fine grids for $kz \gg 1$, we prefer to use a well-conditioned more general BVP solver. 

Now, let us expand all functions in truncated Chebyshev series
\begin{equation}
\label{eq:chebexp}
 y(z) = \sum_{n=0}^{N-1} \four{y}_n, \qquad 
 y^\prime(z) = \sum_{n=0}^{N-1} \four{y}^\prime_n T_n(z),\qquad
y^{\prime \prime} (z) = \sum_{n=0}^{N-1} \four{y}^{\prime \prime}_n T_n(z), \qquad 
f(z) = \sum_{n=0}^{N-1} \four{f}_n T_n(z). 
\end{equation}
To obtain the coefficients $\four{f}_n$, we can evaluate $f$ on a Chebyshev grid with $N$ points, take the periodic extension of $f(z)$, and do a complex FFT  \cite{trefethen2000spectral}. It can be shown, using the indefinite integrals of Chebyshev polynomials, 
\begin{gather}
\int T_0(x) \, dx = T_1(x) + C, \\[2 pt]
\int T_n(x) \, dx = \frac{1}{2}\left(\frac{T_{n+1}(x)}{n+1}-\frac{T_{n-1}(x)}{n-1}\right)+C, \qquad n \geq 2, 
\end{gather}
where $C$ is constant, that
\begin{gather}
\nonumber
\four{y}'_1= \frac{1}{2}\left(2\four{y}''_0-\four{y}''_2\right)\\[2 pt]
\label{eq:ds}
\four{y}'_n = \frac{1}{2n}\left(\four{y}''_{n-1}-\four{y}''_{n+1}\right), \quad n \geq 2. 
\end{gather}

Computing the coefficients of the function $y(z)$, we have the relationships
\begin{gather}
\nonumber
\four{y}_1 = \frac{1}{2}\left(2\four{y}'_0-\four{y}'_2\right)=\four{y}'_0 -\frac{1}{8}\left(\four{y}''_1-\four{y}''_3\right),\\[2 pt]
\nonumber
\four{y}_2 = \frac{1}{4}\left(\four{y}'_1-\four{y}'_3\right)=\frac{1}{4}\left(\frac{1}{2}\left(2\four{y}''_0-\four{y}''_2\right)-\frac{1}{6}\left(\four{y}''_2-\four{y}''_4\right) \right),\\[2 pt]
\label{eq:cs}
\four{y}_n = \frac{1}{2n}\left(\four{y}'_{n-1}-\four{y}'_{n+1}\right)=\frac{1}{2n}\left(\frac{1}{2n-2}\left(\four{y}''_{n-2}-\four{y}''_n\right)-\frac{1}{2n+2}\left(\four{y}''_n-\four{y}''_{n+2}\right) \right), \quad n \geq 3.
\end{gather}
Note that in some presentations, including \cite{greengard1991spectral},
the first Chebyshev coefficient $\four{y}_0$ is defined with a 1/2 in front of it, so that the exceptions in\ \eqref{eq:cs} for $\four{y}_1$ and $\four{y}_2$ can be written as part of the general case. In either case, the formulation gives two free parameters $\four{y}_0$ and $\four{y}'_0$, which are obtained using the boundary conditions. To match the number of unknowns (coefficients) and equations, we assume $\four{y}''_n = \four{y}'_n = 0$ for $n > N-1$ when calculating $\four{y}_n$ using\ \eqref{eq:cs}.  

We can reformulate the boundary value problem using the Chebyshev series representations\ \eqref{eq:chebexp} as
\begin{gather}
\label{eq:bvpn0}
\sum_{n=0}^{N-1} (\four{y}''_n-k^2\four{y}_n)T_n(z) = \sum_{n=0}^{N-1}\four{f}_n T_n(z).
\end{gather}
Matching modes gives a system of equations for the Chebyshev coefficients
\begin{gather}
\label{eq:modes}
\four{y}''_n-k^2\four{y}_n=\four{f}_n \qquad n=0, \dots N-1, \\[2 pt]
\label{eq:BCs}
\sum_{n=0}^{N-1} (\four{y}'_n+k\four{y}_n) = \alpha, \qquad \sum_{n=0}^{N-1} (\four{y}'_n-k\four{y}_n) (-1)^n = \beta.
\end{gather}
with $N+2$ equations and $N+2$ unknowns $\four{y}_0, \four{y}'_0, \four{y}''_0, \dots \four{y}''_{N-1}$. We solve the algebraic system of equations\ \eqref{eq:modes} and \eqref{eq:BCs} for the second derivative coefficients $\four{y}''_0, \dots \four{y}''_{N-1}$ and integration constants $\four{y}_0$ and $\four{y}'_0$, then determine $\four{y}_n$ by applying the double integral operation on the coefficients\ \eqref{eq:cs}.

If the domain is $z \in [0,H]$, then the equations are
 \begin{gather}
\label{eq:BVPH}
\four{y}''_n-k^2\frac{H^2}{4}\four{y}_n=\four{f}_n \qquad n=0, \dots N-1, \\[2 pt]
\label{eq:BVPBCH}
\sum_{n=0}^{N-1} \left(\frac{H}{2}\four{y}'_n+k\frac{H^2}{4}\four{y}_n\right) = \alpha, \qquad \sum_{n=0}^{N-1} \left(\frac{H}{2}\four{y}'_n-k\frac{H^2}{4}\four{y}_n\right) (-1)^n = \beta.
\end{gather}
For the $k=0$ mode, the system reduces to
 \begin{gather}
\sum_{n=0}^{N-1} \four{y}''_n T_n\left(\frac{2z}{H}-1\right) = \sum_{n=0}^{N-1}\four{f}_n T_n\left(\frac{2z}{H}-1\right),
\end{gather}
which gives the trivial system of equations
\begin{gather}
\label{eq:modes0}
\four{y}''_n=\four{f}_n \qquad n=0, \dots N-1.
\end{gather}
Our solver uses homogeneous boundary conditions for the $k=0$ mode, 
\begin{gather}
\label{eq:BCs0}
\sum_{n=0}^{N-1} \four{y}_n = 0, \qquad \sum_{n=0}^{N-1} \four{y}_n (-1)^n = 0. 
\end{gather}

We use a Schur complement approach to solve the algebraic system of equations\ \eqref{eq:BVPH} and\ \eqref{eq:BVPBCH}. We can write the system in block form as 
\begin{equation}
\label{eq:blocksys}
\begin{pmatrix} \bm{A} & \bm{B} \\[2 pt] \bm{C} & \bm{D} \end{pmatrix}
\begin{pmatrix} \four{\bm{y}}''\\ \four{y}_0 \\ \four{y}'_0 \\ \end{pmatrix}
= \begin{pmatrix}  \four{\bm{f}}\\ \alpha \\ \beta\\ \end{pmatrix}. 
\end{equation}
Here $\bm{B}$ is $N \times 2$, $\bm{C}$ is $2 \times N$, $\bm{D}$ is $2 \times 2$, and $\bm{A}$ is an $N \times N$ \textit{pentadiagonal} matrix with only three nonzero diagonals, which we first (pre)factorize using a fast algorithm for such matrices (\cite{karawia2010two}, specialized to the case of only three nonzero diagonals). Our Schur complement approach is then to solve the $2 \times 2$ system
\begin{equation}
(\bm{C}\bm{A}^{-1}\bm{B}-\bm{D})\begin{pmatrix}  \four{y}_0 \\ \four{y}'_0\\ \end{pmatrix} = \bm{C}\bm{A}^{-1}\bm{f}-\begin{pmatrix}  \alpha \\ \beta\\ \end{pmatrix}
\end{equation}
for $\four{y}_0$ and $\four{y}'_0$. We then obtain the coefficients $\four{\bm{y}}''=(\four{y}''_0, \dots \four{y}''_{N-1})$ from back-substitution
\begin{equation}
\four{\bm{y}}''=\bm{A}^{-1}\left(\bm{f}-\bm{B}\begin{pmatrix}   \four{y}_0 \\ \four{y}'_0\\ \end{pmatrix}\right).
\end{equation}

\bibliographystyle{unsrt}

\bibliography{DoublyPeriodicBib,References}

\end{document}